\documentclass{article}
\usepackage{latexsym}
\usepackage{amsfonts}
\usepackage{graphicx}
\usepackage{verbatim}
\usepackage{subcaption}
\usepackage{color}

\newtheorem{theorem}{Theorem}[section]
\newtheorem{lemma}{Lemma}[section]

\newtheorem{definition}{Definition}[section]
\newtheorem{prop}{Proposition}[section]

\newtheorem{uda}{Example}[section]
\newcommand{\qed}{\hfill$\Box$\par\medskip}

\newenvironment{Proof}{\noindent{\sc Proof.}}{\qed}

\def\bhag#1{\noindent
\setcounter{equation}{0}
\section{#1}
}
\def\bfgk#1{{{#1}\kern-5.5pt{#1}}}

\def\RR{{\mathbb R}}
\def\CC{{\mathbb C}}
\def\ZZ{{\mathbb Z}}

\def\PPI{{{\rm I}\kern-1pt\Pi}}
\def\SS{{\mathbb S}}

\def\TT{\mathbb T}

\def\b #1;{{\bf #1}}
\def\x{{\bf x}}
\def\k{{\bf k}}

\def\e{\epsilon}

\def\C{{\mathcal C}}

\def\GG{\mathbb{G}}
\def\YY{\mathbb{Y}}

\def\ip#1#2{{\langle {#1}, {#2}\rangle}}
\def\derf#1#2{{#1}^{(#2)}}

\def\esssup{\mathop{\hbox{{\rm ess sup}}}}

\def\be{\begin{equation}}
\def\ee{\end{equation}}
\def\bea{\begin{eqnarray}}
\def\eea{\end{eqnarray}}
\def\eref#1{(\ref{#1})}
\def\disp{\displaystyle}

\def\donchitre#1#2{\vskip 6.5cm\noindent
\parbox[t]{1in}{\special{eps:#1.eps x=6.5cm y=5.5cm}}
\hbox to 7cm{}\parbox[t]{0.0cm}{\special{eps:#2.eps x=6.5cm y=5.5cm}}}

\def\tn{|\!|\!|}
\def\XX{{\mathbb X}}
\def\BB{{\mathbb B}}
\def\span{\mbox{{\rm span }}}

\title{Minimum Sobolev norm interpolation of derivative data}
\author{ S.~Chandrasekaran
    \thanks{Department of Electrical and
    Computer Engineering, University of California, Santa Barbara,
    Santa Barbara, CA 93106. The research of this author was
    supported, in part, by grants CCF-0515320 and CCF-0830604 from the
    NSF.}
    \and C.~H.~Gorman
    \thanks{Department of Mathematics, University of California, Santa Barbara,
    CA 93106.}
    \and H.~N.~Mhaskar
    \thanks{ Institute of Mathematical Sciences,
    Claremont Graduate University, Claremont, CA 91711
    \textsf{email:} hrushikesh.mhaskar@cgu.edu. The research of this author is supported in part by ARO Grant W911NF-15-1-0385.}
}

\date{}
\begin{document}
\maketitle

\begin{abstract}
We study the problem of reconstructing a function on a manifold satisfying some mild conditions, given data on the values and some derivatives of the function at arbitrary points on the manifold. 
While the problem of finding a polynomial of two variables with total degree $\le n$ given 
the values of the polynomial and some of its derivatives at exactly the same number of points as the dimension of the polynomial space is sometimes impossible, 
we show that such a problem always has a solution in a very general situation if the degree of the polynomials is sufficiently large. We give estimates on how large the degree should be, 
and give explicit constructions for such a polynomial 
even in a far more general case. As the number of sampling points at which the data is available increases,
our polynomials converge to the target function on the set where the sampling points are dense. Numerical examples in single and double precision
  show that this method is stable and of high-order.
\end{abstract}

\bhag{Introduction}
The subject of Lagrange interpolation of univariate functions is a very old one. Thus, one starts with an (infinite) \emph{interpolation matrix} $X$ whose $n$-th column consists of $n$ real numbers $x_{k,n}$, $k=1,\cdots, n$. 
It is well known that for any matrix $Z$ whose $n$-th column consists of $n$ real numbers $z_{k,n}$, $k=1,\cdots,n$, there exists a sequence of polynomials $L_n(X,Z)$ of degree $\le n-1$ such that $L_n(X,Z;x_{k,n})=z_{k,n}$, $k=1,\cdots,n$. In the case when each $z_{k,n}=f(x_{k,n})$ for some continuous function $f$, then it is customary to denote $L_n(X,Z)$ by $L_n(X, f)$.
It is also well known that for any $X\subset [-1,1]$, there exists a continuous function $f :[-1,1]\to\RR$ such that the sequence $L_n(X,f)$ does not converge in the uniform norm~\cite{natanson1982constructive}.
This situation is similar to the theory of trigonometric Fourier series, where the Fourier projections of a function do not always converge to the function in the uniform norm, but one can construct summability operators to obtain convergence~\cite{zygmund}.

Several interpolatory analogues of such summability operators are studied in the literature. 
For example,
if $x_{k,n}=\cos((2k-1)\pi/(2n))$, $k=1,\cdots,n$, $f:[-1,1]\to \RR$ is continuous, and one constructs a sequence of polynomials $F_n(f)$ of degree $\le 2n-1$ such that $F_n(x_{k,n})=f(x_{k,n})$, and $F_n'(x_{k,n})=0$, then the sequence $F_n(f)\to f$ uniformly on 
$[-1,1]$~\cite{natanson1982constructive}.

In 1906, Birkhoff initiated a study of interpolation in a more general setting, known now as Birkhoff interpolation~\cite{Birkhoff1906}. 
For each column of the matrix $X$, one considers a  \emph{incidence matrix} $E$ whose entries are in $\{0,1\}$. If the number of $1$'s in $E$ is $N$, one seeks a polynomial $B_N$ of degree $\le N-1$ such that $B_N^{(j)}(x_{k,n})=y_{j;k,n}$ if the $(k,j)$-th entry in $E$ is $1$. 
Clearly, such a polynomial may or may not exist. The conditions under which it exists and is unique is the topic of a great deal of research~\cite{GLorentz1983}.
As expected, the problems are much harder in the multivariate
setting~\cite{rudy,Lorentz_Review}.
For example, it is not always possible to find a bivariate polynomial $P$ of total degree $2$ (with $6$ parameters) such that $P$ and its derivatives up to order $2$ take given values (also $6$ conditions).  

A great deal of research on this subject is focused on forming a ``square'' system and
determining whether or not the system is solvable for most point
distributions. Since there is not always an obvious choice of
interpolation space, the space can be constructed in multiple
ways. For example, monomial
bases~\cite{CHAI20113207,LEI20111656,CUI2015162}, Newton-type
bases~\cite{WANG2009466,CHAI20113207} and cardinal
bases~\cite{ALLASIA20129248} have been studied by various authors.
Obtaining a unique solution, though, does not ensure convergence to the
underlying function when the data become dense. Although some of the
previously cited articles have numerical examples, not all of them
computed the maximum difference between the test function and its
approximation in a suitable domain, and we do not know of any general
provably convergent technique for the Birkhoff problem on scattered
data points in arbitrary domains.

If we do not require that the dimension of the polynomial space match exactly the number of $1$'s 
in the incidence matrix, then it is possible to guarantee not just existence, but also give explicit algorithms and prove the convergence of the resulting polynomials. 
In the case of interpolation based on
the values of the polynomials alone, this has been observed in a series 
of papers~\cite{bdint,Szabados1990,approxint2002}.
The purpose of this paper is to generalize these results for Birkhoff-like interpolation for the so called diffusion polynomials.

One of our motivations is to study numerical solutions of a system of linear partial differential equations. Given differential operators $L_k$ on a manifold $\XX$ (and its boundary), the collocation method involves finding a ``polynomial'' (i.e., an element of a suitably chosen finite dimensional space)
$P$ for which the values of $L_k(P)$ are known at some grid points. We view this question as a generalized Birkhoff interpolation problem, except that we do not require the dimension of the space to 
be exactly equal to the conditions.
On general manifolds, one does not always have standard grids such as equidistant grid on a Euclidean space. 
Therefore mesh-free methods require a solution of such interpolation problems on scattered data; i.e., when one cannot prescribe the location of the grid points in advance. 
We show the feasibility of the solution of such interpolation problems provided the dimension of the space is sufficiently high, proportional to the minimal separation among the grid points. 
We will prove that such solutions can be constructed as minimizers of an optimization problem, and prove that these solutions will converge to the target function at limit points of the grid points. 
An application of these ideas is already demonstrated 
in~\cite{chandra2015minimum}.
We will
give several numerical examples to illustrate the concepts introduced in this paper in the context of the sphere.

The main results are presented in Section~\ref{mainsect}, while the
overall assumptions are discussed in Section~\ref{asssect}.
Preparatory results are developed in
Sections~\ref{pfprepsect}--\ref{approxsect} and are used to prove the
main theorems in Section~\ref{pfsect}.  Finally, numerical simulations
are given in Section~\ref{num_results}.  We show that the standard
divergence phenomenon is overcome in both one and two dimensions. We
also numerically demonstrate the high-order convergence of our method.
Furthermore we present results using both \textsl{single} and double
precision to show that the high ill-conditioning associated with high
order methods can be successfully overcome using our numerical
techniques.

\bhag{Main results}\label{mainsect}

We wish to study the problem of interpolation of derivative information
in a somewhat abstract manner, to accommodate several examples.
The most elementary among these is interpolation by multivariate
trigonometric polynomials. Other examples include the problem
of interpolation by spherical polynomials, or by linear combinations
of the eigenfunctions of some elliptic differential operator
on a smooth manifold. 


Let $\XX$ be a  metric measure space with a
probability measure $\mu^*$ and metric $\rho$. In this paper, a measure will mean a complete, sigma-finite, Borel measure, signed or positive.  For a measure $\nu$, $|\nu|$ denotes its total variation measure. If $\nu$ is a measure, and 
$f :\XX\to\RR$ is $\nu$-measurable, we define
$$
\|f\|_{\nu;p}=\left\{\begin{array}{ll}
\disp\left(\int_\XX |f(x)|^pd|\nu|(x)\right)^{1/p}, &\mbox{ if $1\le p<\infty$,}\\
\disp|\nu|-\esssup_{x\in\XX}|f(x)|, &\mbox{ if $p=\infty$.}
\end{array}\right.
$$
The symbol $L^p(\nu)$ will denote the space of all $f$ such that
$\|f\|_{\nu;p}<\infty$, where two functions are considered equal if they
are equal $|\nu|$--almost everywhere. If $\nu=\mu^*$, we will often omit the mention of the measure from the notation, if we feel that this should not cause any confusion; e.g., $\|f\|_p=\|f\|_{\mu^*;p}$, $L^p=L^p(\mu^*)$. The space of uniformly continuous and bounded functions on $\XX$ will be denoted by $UBC$. For $1\le p\le\infty$, we define
the dual exponent $p'$ by $1/p+1/p'=1$ as usual.

Let $\{\lambda_k\}$ be increasing sequence of nonnegative numbers
with $\lambda_0=0$, and $\{\phi_k\}$ be an orthonormal set (of real
valued functions) in $L^2\cap UBC$.  We define the space
$$
\Pi_n = \span\{\phi_k : \lambda_k<n\}, \qquad n>0,\ \Pi_n=\{0\}, \qquad n\le 0.
$$
We will write
$$
\Pi_\infty =\bigcup_{n\ge 0}\Pi_n.
$$
The elements of $\Pi_\infty$ have been referred to as diffusion
polynomials in \cite{mauropap}, and we will use the same terminology.
If $P\in\Pi_n$, we will refer to $n$ as the degree of $P$, or more
precisely that $P$ is of degree $<n$. The $L^p$ closure of the
$\Pi_\infty$ will be denoted by $X^p$. 

\begin{uda}\label{trigexam1}
{\rm
\textbf{The trigonometric case}
Let $q\ge 1$ be an integer. The space $\XX$ is the torus
(quotient group) $\TT^q=\RR^q/(2\pi\ZZ^q)$, $\rho$ is the
arc--length, and $\mu^*$ is the normalized Lebesgue measure.
We take $\{\lambda_k\}$ to be an enumeration of multi--integers in $\ZZ_+^q$
such that $\mathbf{m}$ comes before $\mathbf{j}$ if either
$|\mathbf{m}|_2\le |\mathbf{j}|_2$ and if
$|\mathbf{m}|_2=|\mathbf{j}|_2$, then $\mathbf{m}$ comes before
$\mathbf{j}$ in the alphabetical ordering.
If the multi--integer corresponding to $\lambda_k$ is $\k$,
the corresponding $\phi_k$'s are $\cos(\k\cdot\circ), \sin(\k\cdot\circ)$. In this case,
$X^p=L^p$ if $1\le p<\infty$, and $X^\infty$ is the set of all
 continuous functions on $\TT^q$.
}
\end{uda}
\begin{uda}\label{manifoldexam1}
{\rm \textbf{The manifold case.}
Let $\XX$ be a compact Riemannian manifold, $\rho$ be the geodesic
distance, $\mu^*$ be the Riemann measure. We take $\lambda_k$ to be the
square roots of eigenvalues of a self--adjoint, second order,
regular elliptic differential operator on $\XX$, and $\phi_k$
to be the corresponding eigenfunction. The precise nature of the
space $X^p$ will depend upon the manifold and the elliptic operator. 

In particular, if $\XX$ is the $2$ dimensional Euclidean sphere $\SS^2$, 
we may choose each $\phi_k$ to be one of the orthonormalized
spherical harmonics, and $\lambda_k$ to be the degree of this polynomial.
}
\end{uda}

 In the absence of any concrete structure, we will need to make
several assumptions on the system
$\Xi=(\XX,\rho,\mu^*, \{\lambda_k\}, \{\phi_k\})$.
These are formulated precisely in Section~\ref{asssect}.
For the clarity of exposition, we will now assume that these
are all satisfied. In particular, the symbol $q$ used in the
following discussion is defined in \eref{ballmeasurecond}.
 
 Next, we define the smoothness classes needed in order to state
our theorems. In the trigonometric case, the Sobolev space
$W_{p,\beta}$ is defined as the space of all functions
$f:\TT^q\to\CC$ for which
$(|\k|_2^2+1)^{\beta/2}\hat{f}(\k)=\widehat{f^{(\beta)}}(\k)$
for some $f^{(\beta)}\in L^p(\TT^q)$. In our abstract case,
the role of $|\k|_2$ is played by $\lambda_k$.
However, we would like to introduce a greater flexibility in
the definition of the Sobolev class. 

If $f\in L^1$, we  define
\be\label{hatfdef}
\hat f(k):=\int_\XX f(y){\phi_k(y)}d\mu^*(y), \qquad k=0,1,\cdots.
\ee
To include such multipliers as $(\lambda_k^2+1)^{\beta/2}$ or
$(\lambda_k+1)^\beta$, we use a mask of type $\beta$, defined below in
Definition~\ref{eigkerndef}.
In the sequel, $S$ will be a fixed integer.

\begin{definition}\label{eigkerndef}
Let $\beta\in\RR$. A function $b:\RR\to\RR$ will be called a mask
of type $\beta$ if $b$ is an even, $S+1$ times continuously
differentiable function (for some integer $S>q$) such that for $t>0$,
$b(t)=(1+t)^{-\beta}F_b(\log t)$ for some $F_b:\RR\to\RR$ such that
$|\derf{F_b}{k}(t)|\le c(b)$, $t\in\RR$, $k=0,1,\cdots,S+1$,
and $F_b(t)\ge c_1(b)$, $t\in\RR$. 
\end{definition}
 
If $\beta\in\RR$, $b$ is a mask of type $\beta$, and $f\in L^p$,
we say that a function $f\in W_{p,b}$ if there exists $f^{(b)}\in L^p$
such that
\be\label{wbspacedef}
b(\lambda_k)\widehat{f^{(b)}}(k)=\hat{f}(k), \qquad k=0,1,\cdots.
\ee
We will write
\be\label{wbnormdef}
\|f\|_{W_{p,b}}=\|f^{(b)}\|_p.
\ee
In \cite{eignet}, we have shown that if $\beta>q/p$,
and $b$ is a mask of type $\beta$, then for every $y\in\XX$,
there exists $\psi_y:=G(b;\circ,y)\in X^{p'}$ such that
$\ip{\psi_y}{\phi_k} =b(\lambda_k)\phi_k(y)$, $k=0,1,\cdots$. Moreover,
\be\label{glpuniform}
\sup_{y\in\XX}\|G(b;\circ,y)\|_{p'} \le c.
\ee
 It is then easy to check by comparison of coefficients that 
if $f\in W_{p,b}$ then
for almost all $x\in\XX$:
\be\label{reproducekerndef}
f(x) =\int_\XX G(b;x,y)\derf{f}{b}(y)d\mu^*(y).
\ee

As a prelude to our main theorem, we formulate first a
Golomb--Weinberger--type theorem in our general setting.
The theorem gives an explicit expression for a solution of the
interpolation problem:\\

\noindent
\emph{Given linear operators $L_k$, $k=1,\cdots, R$,
each defined on $W_{2,b}$, such that point evaluations are
well--defined on the range of each $L_k$, points
$y_\ell\in\XX$, $\ell=1,\cdots,M$, and data $f_{k,\ell} \in\CC$,
find a function $g\in W_{2,b}$ such that
$L_k(g)(y_\ell)=f_{k,\ell}$, $k=1,\cdots, R$, $\ell=1,\cdots, M$.}\\

The Golomb--Weinberger--type theorem gives a solution as a
linear combination of the kernel defined in \eref{kerndef} below.
Let $b$ be a mask of type $\beta>q/2$. Then $b^2$ is a mask of type
$2\beta>q$. We  define
\be\label{kerndef}
\GG(x,y)=\sum_{j=0}^\infty b(\lambda_j)^2\phi_j(x){\phi_j(y)}.
\ee
Clearly, $\GG$, treated either as a function of $x$ or $y$,
is a function in $W_{2,b}$.

\begin{theorem}\label{golombtheo}
Let $R, M\ge 1$ be integers, each $L_k$, $k=1,\cdots,R$, 
be a linear operator  defined on $W_{2,b}$ such that point evaluations
are well defined on the range of $L_k$. Let $y_\ell\in\XX$,
$\ell=1,\cdots,M$, and
$\{f_{k,\ell}\}_{k=1,\cdots,R,\ \ell=1,\cdots, M} \subset \CC$.
We assume that there exist
$\phi_{k,j}\in W_{2,b}$, $k=1,\cdots,R$, $j=1,\cdots, M$,
such that  the following condition holds:
\textit{ For $i=1,\cdots,R$, $\ell=1,\cdots,M$, 
\be\label{dualfndef}
L_i(\phi_{k,j})(y_\ell)=\left\{\begin{array}{ll}
1, &\mbox{ if $i=k$, $j=\ell$,}\\
0, & \mbox{ otherwise. }
\end{array}\right.
\ee
}
Then the problem
\be\label{golomb}
\mbox{ minimize }\|\derf{g}{b}\|_2 \mbox{ subject to }
L_k(g)(y_\ell)=f_{k,\ell},
  \qquad k=1,\cdots,R,\ \ell=1,\cdots,M, 
\ee
has a solution of the form
\be\label{kernsol}
P(x)=\sum_{k=1}^R\sum_{j=1}^{M} a_{k,j}L_{k,1}\GG(y_j,x), \qquad x\in\XX,
\ee
where, the expression $ L_{k,1}\GG(y_j,x)$ means that the operator
$L_k$ is applied to the $1$--st  variable in $\GG$,
and the resulting function is evaluated at $(y_j,x)$.
\end{theorem}

In the trigonometric case, when $b(t)=(t^2+1)^{-\beta/2}$,
the kernel $\GG$ takes the form
$$
\sum_{\k\in\ZZ}(|\k|_2^2+1)^{-\beta}\exp(i\k\cdot\x).
$$
A straightforward computation of this series may be slow for
small values of $\beta$, and in any case, introduces a truncation error.
Therefore, we are interested in solving the interpolation problem
directly using the diffusion polynomials (trigonometric polynomials
in the trigonometric case).
Toward this goal, we first introduce some further terminology.

\noindent
\textbf{Constant convention}:\\
\emph{In the sequel, $c, c_1,\cdots$ will denote generic positive
constants independent of any obvious variables such as the degree $n$,
or the target function, etc. Their values may be different at
different occurences, even within the same formula.
The symbol $A\sim B$ means that $c_{1}A\le B\le c_2A$. }\\

We will consider an \emph{interpolation matrix}, by which we mean a
sequence $\{Y_n\}_{n=1}^\infty$ of subsets 
$$
Y_n=\{y_{j,n}\}_{j=1}^{M_n}, \qquad n=1, 2, \cdots.
$$
This is not a matrix in the usual sense of the word, but the
terminology captures the spirit of a similar notion in the theory of
classical polynomial  interpolation. 
For any subset $\C\subseteq\XX$, we define its minimal separation by
\be\label{minsepdef}
\eta(\C)=\min_{x, y\in\C,\ x\not=y}\rho(x,y).
\ee
We will write
\be\label{minsepspecific}
\eta_n=\eta(Y_n), \qquad n=1,\cdots.
\ee

Our main theorem is the following.

\begin{theorem}\label{feasibletheo}
We assume each of the  conditions listed in Section~\ref{asssect},
where some of the notation used here is explained.
Let $1\le p \le \infty$, $\beta>\max_{1\le k\le R}q_k +q/p$, $b$
be a mask of type $\beta$, $f\in W_{p,b}$.
Then there exists an integer $N^*$ with $N^*\sim \eta_n^{-1}$
  and a mapping ${\bf P}^*=\mathbf{P}^* :W_{p,b}\to \Pi_{N^*}$
such that for every
  $f\in W_{p,b}$, 
\be\label{sobinterp}
 L_k({\bf P}^*(f))(y_{j,n})=L_k(f)(y_{j,n}),
  \qquad j=1,\cdots,M_n, \ k=1,\cdots,R,
\ee 
and 
\be\label{sobapprox}
 \|f-{\bf  P}^*(f)\|_{W_{p,b}}\le c \inf\{ \|f-T\|_{W_{p,b}}\ :\ T\in
  \Pi_{N^*}\}.  
\ee 
In particular, there exists
$\mathbf{P}=\mathbf{P}_n :W_{p,b}\to \Pi_{N^*}$ such that
\be\label{sobminnorm}
\|\mathbf{P}(f)\|_{W_{p.b}}=\min\{\|P\|_{W_{p,b}} :
L_k(P)(y_{\ell,n})=L_k(f)(y_{\ell,n}),\ \ell=1,\cdots, M_n,
\ k = 1,\cdots, R\}.
\ee
\end{theorem}

We observe that the sets $Y_n$ do not necessarily become dense in
$\XX$ as $n\to\infty$.
Theorem~\ref{feasibletheo} helps us to find an interpolatory
diffusion polynomial whose $W_{p,b}$ norm is under control on $\XX$.
Therefore, we do not expect that the sequence $\mathbf{P}_n(f)$
to converge to $f$ on $\XX$. However, the sequence converges at
limit points of $Y_n$'s.

\begin{theorem}\label{feasimpconvtheo}
With the set up as in Theorem~\ref{feasibletheo}, if $x_0\in\XX$
is a limit point of the family $\{Y_n\}$, then $f(x_0)$ is a
limit point of $\{\mathbf{P}_n(f)(y) : y\in Y_n, \ n=1,2,\cdots\}$.
\end{theorem} 

Theorem~\ref{feasimpconvtheo} follows from a much general principle
``feasibility implies convergence'', which is formulated more
precisely in Theorem~\ref{convergetheo} below.

\bhag{Assumptions}\label{asssect}

\subsection{The space}\label{spaceasssect}
Let $\XX$ be a non-empty set, $\rho$ be a metric defined on $\XX$, and $\mu^*$ be a complete, positive, Borel measure with $\mu^*(\XX)=1$.
We fix a non-decreasing sequence $\{\lambda_k\}_{k=0}^\infty$ of nonnegative numbers such that $\lambda_0=0$, and $\lambda_k\uparrow \infty$ as $k\to\infty$. Also, we fix a system of continuous, bounded, and integrable functions $\{\phi_k\}_{k=0}^\infty$, orthonormal with respect to $\mu^*$; namely, for all nonnegative integers $j, k$,
\begin{equation}\label{orthonormality}
\int_\mathbb{X} \phi_k(x)\phi_j(x)d\mu^*(x) =\left\{\begin{array}{ll}
1, &\mbox{if $j=k$,}\\
0, &\mbox{otherwise.}
\end{array}\right.
\end{equation}
We will assume  that $\phi_0(x)=1$ for all $x\in\mathbb{X}$. 

In our context, the role of polynomials will be played by diffusion polynomials, which are finite linear combinations of $\{\phi_j\}$. In particular, an element of
$$
\Pi_n :=\mathsf{span}\{\phi_j : \lambda_j <n\}
$$
will be called a diffusion polynomial of degree $<n$.

We will formulate our assumptions in terms of a formal heat kernel. The \emph{heat kernel} on $\mathbb{X}$ is defined formally by
\begin{equation}\label{heatkerndef}
K_t(x,y)=\sum_{k=0}^\infty \exp(-\lambda_k^2t)\phi_k(x){\phi_k(y)}, \qquad x, y\in \mathbb{X}, \ t>0.
\end{equation}
Although $K_t$ satisfies the semigroup property, and in light of the fact that $\lambda_0=0$, $\phi_0(x)\equiv 1$, we have formally
\begin{equation}\label{heatkernint}
\int_\mathbb{X} K_t(x,y)d\mu^*(y) =1, \qquad x\in\mathbb{X},
\end{equation}
yet $K_t$ may  not be the heat kernel in the classical sense. In particular, we need not assume $K_t$ to be nonnegative.

\begin{definition}\label{ddrdef}
The system $\Xi=(\mathbb{X},\rho,\mu^*,\{\lambda_k\}_{k=0}^\infty,\{\phi_k\}_{k=0}^\infty)$ is called a \textbf{data-defined space} if each of the following conditions are satisfied.
\begin{enumerate}
\item For each $x\in\mathbb{X}$ and $r>0$, the ball $\mathbb{B}(x,r)$ is compact.
\item There exist $q>0$ and $\kappa_2>0$ such that the following power growth bound condition holds:
\begin{equation}\label{ballmeasurecond}
\mu^*(\mathbb{B}(x,r))=\mu^*\left(\{y\in\mathbb{X} : \rho(x,y)<r\}\right) \le \kappa_2r^q, \qquad x\in\mathbb{X}, \ r>0.
\end{equation}
\item The series defining $K_t(x,y)$ converges for every $t\in (0,1$ and $x,y\in\mathbb{X}$. Further, with $q$ as above, there exist $\kappa_3, \kappa_4>0$ such that the following Gaussian upper bound holds:
\begin{equation}\label{gaussianbd}
|K_t(x,y)|\le \kappa_3t^{-q/2}\exp\left(-\kappa_4\frac{\rho(x,y)^2}{t}\right), \qquad x, y\in \mathbb{X},\ 0<t\le 1.
\end{equation} 
\end{enumerate}
\end{definition}

There is a great deal of discussion in the literature on the validity of the  conditions in the above definition and their relationship with many other objects related to the quasi--metric space in question, (cf. for example,  \cite{davies1997,  grigor1999estimates, grigor2006heat, grigorlyan2heat}). 
In particular, it is shown in \cite[Section~5.5]{davies1997} that all the conditions defining a data-defined space are satisfied in the case of any complete, connected Riemannian manifold with non--negative Ricci curvature. 
It is shown in  \cite{kordyukov1991p}  that  our assumption on the heat kernel is valid in the case when $\mathbb{X}$ is a complete Riemannian manifold with bounded geometry, and $\{-\lambda_j^2\}$, respectively $\{\phi_j\}$, are eigenvalues, respectively eigenfunctions, for a uniformly elliptic second order differential operator satisfying certain technical conditions.

The bounds on the heat kernel are closely connected with the measures of the balls $\mathbb{B}(x,r)$. For example, using \eref{gaussianbd}, Proposition~\ref{criticalprop} below, and the fact that
$$
\int_{\mathbb{X}} |K_t(x,y)|d\mu^*(y)\ge \int_{\mathbb{X}} K_t(x,y)d\mu^*(y)=1, \qquad x\in \mathbb{X},
$$
it is not difficult to deduce as in \cite{grigorlyan2heat} that
\begin{equation}\label{ballmeasurelowbd}
\mu^*(\mathbb{B}(x,r))\ge cr^q, \qquad 0<r\le 1.
\end{equation}
In many of the examples cited above, the kernel $K_t$ also satisfies a lower bound to match the upper bound in \eref{gaussianbd}. In this case, Grigory\'an \cite{grigorlyan2heat} has also shown that \eref{ballmeasurecond} is satisfied for $0<r<1$. 

We remark that the estimates  \eref{ballmeasurecond} and \eref{ballmeasurelowbd} together imply that $\mu^*$ satisfies the homogeneity condition
\begin{equation}\label{doublingcond}
\mu^*(\mathbb{B}(x,R))\le c_1(R/r)^q\mu^*(\mathbb{B}(x,r)), \qquad x\in\mathbb{X},\  r\in (0,1],\ R>0,
\end{equation}
where $c_1>0$ is a suitable constant.

\subsection{The operators $L_k$}\label{opasssect}
In this sub--section, we state our assumptions on the linear operators $L_k$.
For a bivariate function $F :\XX\times\XX\to \RR$, we will denote 
$$
L_{k,1}F(x,y):=L_kF(x,y):=(L_k(F(\circ,y))(x),
\qquad L_{k,2}F(x,y):=(L_k(F(x,\circ))(y).
$$
We fix $n$, letting $\eta = \eta_{n}$ be the minimum separation between
the $M = M_{n}$ points in $Y_{n}$.

\begin{enumerate}
\item We assume that each $L_k$ is \emph{closed} linear operator; i.e.,
if $f_m\to f$ in $L^p$, and $L_k(f_m)\to g$ in $L^p$,
then $f$ is in the domain of $L_k$ and $L_k(f)=g$. 
\item We assume  that each $L_k$ is \emph{local};
i.e., if $f(x)=0$ for almost all $x$ in an open subset $U$ of $\XX$,
then $L_k(f)(x)=0$ for almost all $x\in U$. 
\item There exists $q_k\ge 0$ such that 
\be\label{diffgaussianbd}
|L_{k,1}K_t(x,y)| \le ct^{-(q+q_k)/2}\exp\left(-c_1\frac{\rho(x,y)^2}{t}\right), \qquad x, y\in \mathbb{X},\ 0<t\le 1.
\ee
\item For each $k=1,\cdots,R$, $j=1,\cdots, M$, there exists
$\phi_{k,j}\in W_{p,b}$ such that each of the following conditions holds:
\begin{enumerate}
\item For $i=1,\cdots$, 
\be\label{dualityop}
L_i(\phi_{k,j})=0, \qquad i\not=k, \ i=1,\cdots, R.
\ee
\item For $i=1,\cdots$, 
$\ell=1,\cdots,M$, 
\be\label{dualfndefbis}
L_i(\phi_{k,j})(y_\ell)=\left\{\begin{array}{ll}
1, &\mbox{ if $i=k$, $j=\ell$,}\\
0, & \mbox{ otherwise. }
\end{array}\right.
\ee
\item $\phi_{k,j}$, $\phi_{k,j}^{(b)}$ are both supported on a
neighborhood of $y_j$ with diameter $\le \eta/3$.
\item We have
\be\label{dualfnnormest}
\|\phi_{k,j}^{(b)}\|_\infty \le c\eta^{q_k-\beta}.
\ee
\end{enumerate}

\end{enumerate}

We remark that  for any $x\in \XX$, there is at most one $j$ such that
$\rho(x,y_j)\le \eta/3$, and hence, 
\be\label{dualitypt}
\phi_{k,\ell}(x)=0, \qquad \ell\not=j,\ j=1,\cdots,M.
\ee

In the case when $\XX$ is the torus, and each $L_k$ is a mixed
partial derivative $D^{{\bf j}}$ for some multi--integer ${\bf j}$,
the construction of functions $\phi_{k,j}$ is given in \cite{sharmapap}.
The bounds \eref{dualfnnormest} are established there in the case
when $\beta$ is an even integer.
If $\XX$ is a manifold (sphere in particular), and $\eta$ is less
than its inradius, $(\lambda_k^2,\phi_k)$ are the eigenfunctions
of the Laplace--Beltrami operator, and $\beta$ is an even integer,
then the same construction works via exponential coordinates.
A slightly more general scenario involving other second order
partial differential operators holds also in view of our
results in \cite{frankbern}. When the operators $L_k$ are given by
\be\label{pdeopform}
L_k(f,x)=\sum_{|{\bf j}|\le q_k} a_{k,{\bf j}}(x)D^{{\bf j}}f(x),
\ee
then one needs to make some assumptions on the matrices
$(a_{k,{\bf j}}(y_\ell))_{\ell=1,\cdots,M, |{\bf j}|\le q_k}$
so that the constructions given in \cite{sharmapap} can be used
to construct the desired $\phi_{k,j}$. 

\bhag{Preparatory results}\label{pfprepsect}
The proof of Theorem~\ref{feasibletheo} is much more involved
than those of the other theorems. The goal of this section is
to prove a number of auxiliary results which will lead to the
proof of  Theorem~\ref{feasibletheo}.

\subsection{Regular measures}\label{regmeassect}
We start by  introducing the concept of what we have called
$d$--regular measures, and some of their properties.  

\begin{definition}\label{regularitydef}
Let $\nu$ be any  measure on $\XX$ with $|\nu|(\XX)<\infty$, $d>0$.
 We say that $\nu$ is \textrm{$ d$--regular} if
\be\label{regulardef}
\nu(\BB(x,d))\le cd^q, \qquad x\in\XX.
\ee
The infimum of all constants $c$ which work in \eref{regulardef}
will be denoted by $\tn\nu\tn_{R,d}$. 
\end{definition}

For example, \eref{ballmeasurecond} implies that $\mu^*$ is
$d$--regular for every $d>0$, and $\tn\mu^*\tn_{R,d}\le c$
for all $d>0$ with $c$ independent of $d$. In the sequel, when thinking of $\mu^*$ as a regular measure, we will pass to a limit, and use $d=0$.
 In the following lemma, we give another example. 

\begin{lemma}\label{minsepmeaslemma}
 Let $\C=\{y_1,\cdots,y_M\}\subset \XX$, $0<\eta<2$ be the minimal
separation amongst the $y_j$'s (cf. \eref{minsepdef}), $\tau$ be the
measure that associates the mass $\eta^q$ with each $y_j$.
Then $\nu$ is $\eta$--regular, and $\tn\tau\tn_{R,\eta} \le c$,
where $c$ is a constant independent of $\eta$.
\end{lemma}

\begin{Proof} \ 
Let $x_0\in\XX$, and by relabeling if necessary,
let $\C\cap \BB(x_0,\eta)=\{y_1,\cdots,y_J\}$.
Then the caps $\BB(y_j,\eta/2)$ are mutually disjoint,
and their union is contained in $\BB(x_0,3\eta/2)$.
We recall from \eref{ballmeasurelowbd} that
$\eta^q\le c_1\mu^*(\BB(y_j,\eta/2))$ and from
\eref{ballmeasurecond} that $\mu^*(\BB(x_0,3\eta/2))\le c_2\eta^q$.
Therefore, we deduce that
$$
\tau(\BB(x_0,\eta))=J\eta^q \le c_1\sum_{j=1}^J
\mu^*(\BB(y_j,\eta/2))=c_1\mu^*\left(\cup_{j=1}^J \BB(y_j,\eta/2)\right)
\le c_1\mu^*(\BB(x_0,3\eta/2))\le c_1c_2\eta^q.
$$
\end{Proof}

The following proposition
\cite[Theorem~5.5(a), Proposition~5.6]{modlpmz} lists some
equivalent conditions for a measure to be a regular measure.
\begin{prop}\label{mzmeasureprop}
Let $N, d>0$, $\nu$ be a signed or positive finite Borel measure.\\
\noindent
{\rm (a)} If $\nu$ is $d$--regular, then for each $r>0$ and $x\in\XX$,
\be\label{regreconcile}
|\nu|(\BB(x,r))\le c\tn\nu\tn_{R,d}\ \mu^*(\BB(x,r+d))\le 
c_1\tn\nu\tn_{R,d}(r+d)^q.
\ee
Conversely, if for some $A>0$, $|\nu|(\BB(x,r))\le A(r+d)^q$
for each $r>0$ and $x\in\XX$, then $\nu$ is $d$--regular, and
$\tn\nu\tn_{R,d}\le 2^q A$.\\
{\rm (b)} For each $\gamma>1$, 
\be\label{regequiv}
\tn\nu\tn_{R,\gamma d}\le c_1(\gamma+1)^q \tn\nu\tn_{R,d}\le
c_1(\gamma+1)^q\gamma^q\tn\nu\tn_{R,\gamma d},
\ee
where $c_1$ is the constant appearing in \eref{regreconcile}.\\
{\rm (c)} Let $N\ge 1$. If $\nu$ is $1/N$--regular, then $\|P\|_{\nu;p}\le
c_1\tn\nu\tn_{R,1/N}^{1/p}\|P\|_{\mu^*;p}$ for all $P\in \Pi_N$ and
$1\le p<\infty$.
Conversely, if for some $A>0$ and $1\le p<\infty$,
$\|P\|_{\nu;p}\le A^{1/p}\|P\|_{\mu^*;p}$ for all $P\in \Pi_N$,
then $\nu$ is $1/N$--regular, and $\tn\nu\tn_{R,1/N}\le c_2A$.
\end{prop}

Next, we recall a very general proposition \cite[Proposition~6.5]{modlpmz}
helping us to estimate integrals of quantities such as the right
hand side of \eref{opkernlocest}.

\begin{prop}\label{criticalprop}
Let $ d>0$,  and $\nu$ be a $d$--regular measure. 
 If $g_1:[0,\infty)\to [0,\infty)$ is a nonincreasing function,
then for any $N>0$, $r>0$, $x\in\XX$,
\be\label{g1ineq}
N^q\int_{\Delta(x,r)}g_1(N\rho(x,y))d|\nu|(y)\le
\frac{2^{q}(\kappa_1+(d/r)^q)q}{1-2^{-q}}\tn
\nu\tn_{R,d}\int_{rN/2}^\infty g_1(u)u^{q-1}du.
\ee
In particular, if $S>q$, 
\be\label{farptest}
N^q\int_{\Delta(x,r)}\frac{d|\nu|(y)}{\max(1, (N\rho(x,y))^S)}
\le \frac{c_1(c+(d/r)^q)q}{\max(1, (rN/2)^{S-q})}\tn\nu\tn_{R,d},
\ee
and
\be\label{totalnormest}
N^q\int_\XX \frac{d|\nu|(y)}{\max(1, (N\rho(x,y))^S)} \le
c_1(c+(Nd)^q)\tn\nu\tn_{R,d}.
\ee
\end{prop}

\subsection{Localized kernels}\label{lockernsect}

Localized kernels form the main ingredient in our proofs. 
To obtain these kernels, we will use the following Tauberian theorem 
(\cite[Theorem~4.3]{tauberian}) with different choices of the function $H$.

\begin{theorem}\label{maintaubertheo}
Let $\mu^*$ be an extended complex valued measure on $[0,\infty)$, and $\mu^*(\{0\})=0$. We assume that there exist $Q, r>0$, such that each of the following conditions are satisfied.
\begin{enumerate}
\item 
\be\label{muchristbd}
\tn\mu^*\tn_Q:=\sup_{u\in [0,\infty)}\frac{|\mu^*|([0,u))}{(u+2)^Q} <\infty,
\ee
\item There are constants $c, C >0$,  such that
\be\label{muheatgaussbd}
\left|\int_\RR \exp(-u^2t)d\mu^*(u)\right|\le c_1t^{-C}\exp(-r^2/t)\tn\mu^*\tn_Q, \qquad 0<t\le 1.
\ee 
\end{enumerate}
Let $H:[0,\infty)\to\RR$, $S>Q+1$ be an integer, and suppose that there exists a measure $H^{[S]}$ such that
\be\label{Hbvcondnew}
H(u)=\int_0^\infty (v^2-u^2)_+^{S}dH^{[S]}(v), \qquad u\in\RR,
\ee
and
\be\label{Hbvintbdnew}
V_{Q,S}(H)=\max\left(\int_0^\infty (v+2)^Qv^{2S}d|H^{[S]}|(v), \int_0^\infty (v+2)^Qv^Sd|H^{[S]}|(v)\right)<\infty.
\ee
Then for $n\ge 1$,
\be\label{mulockernest}
\left|\int_0^\infty H(u/n)d\mu^*(u)\right| \le c\frac{n^Q}{\max(1, (nr)^S)}V_{Q,S}(H)\tn\mu^*\tn_Q.
\ee
\end{theorem}

We observe that if $H$ is compactly supported, has $S+1$ continuous derivatives,  and is constant in a neigborhood of $0$, then the Taylor formula used with $u\to H(\sqrt{u})$ shows that the representation \eref{Hbvcondnew} holds with $H^{[S]}$ being the $S$-th derivative of $u\to H(\sqrt{u})$, and we have
\be\label{vsqest}
V_{Q,S}(H)\le c\max_{0\le k \le S+1}\max_{u\in\RR}|H^{(k)}(u)|=c|\|H\||_S.
\ee
The following proposition summarizes some general results which we will use in our proofs later. If $\{\psi_j\}$, $\{\tilde{\psi}_j\}$ are sequences of bounded functions on $\XX$, we define formally
\be\label{formalkerndef}
\Phi_n(\{\psi_j\}, \{\tilde{\psi}_j\}, H; x,y)=\sum_{j=0}^\infty H\left(\frac{\lambda_j}{n}\right)\psi_j(x){\tilde{\psi}_j(y)}.
\ee
It is convenient to set 
$$
\Phi_0(\{\psi_j\}, \{\tilde{\psi}_j\}, H; x,y)=\sum_{j: \lambda_j=0} \psi_j(x){\tilde{\psi}_j(y)}.
$$
\begin{prop}\label{acttauberprop}
Let $\{\psi_j\}$, $\{\tilde{\psi}_j\}$ be sequences of bounded functions on $\XX$, $H$ be as in Theorem~\ref{maintaubertheo}, and for $x, y\in\XX$, 
\be\label{genchristbd}
\sum_{j: \lambda_j<u}|\psi_j(x)\tilde{\psi}_j(y)|\le cu^Q, \qquad u\ge 1,
\ee
\be\label{gengaussbd}
\sum_{j=0}^\infty \exp(-\lambda_j^2t)\psi_j(x){\tilde{\psi}_j(y)} \le 
c_1t^{-C}\exp(-\rho(x,y)^2/t), \qquad 0<t\le 1.
\ee
Then
\be\label{genlockernest}
\left|\Phi_n(\{\psi_j\}, \{\tilde{\psi}_j\}, H; x,y)\right|\le c\frac{n^Q}{\max(1, (n\rho(x,y))^S)}V_{Q,S}(H), \qquad n\ge 1.
\ee
Further, if $d>0$ and $\nu$ is a $d$-regular measure, then for $x\in\XX$, 
$r>0$, $n\ge 1$ and $1\le p<\infty$,
\be\label{phinfarest}
\int_{\Delta(x,r)}|\Phi_n(\{\psi_j\}, \{\tilde{\psi}_j\}, H; x,y)|d|\nu|(y)\le 
 n^{Q-q}\frac{c_1(c+(d/r)^q)q}{\max(1, (rN/2)^{S-q})}\tn\nu\tn_{R,d}V_{Q,S}(H),
\ee
\be\label{phinnormest}
\int_\XX |\Phi_n(\{\psi_j\}, \{\tilde{\psi}_j\}, H; x,y)|^pd|\nu|(y)\le c\left(c_1+(nd)^q)\right)n^{Qp-q}\tn\nu\tn_{R,d}V_{Q,S}(H)^p.
\ee
\end{prop}
\begin{Proof}\ 
For each $x, y\in \XX$, we apply Theorem~\ref{maintaubertheo} with the measure
$$
\mu^*_{x,y}([0,u))=\sum_{j: \lambda_j<u}\psi_j(x){\tilde{\psi}_j(y)}.
$$
The estimate \eref{genchristbd} (respectively, \eref{gengaussbd}, \eref{genlockernest}) is equivalent to \eref{muchristbd}  (respectively, \eref{muheatgaussbd}, \eref{mulockernest}). The estimates \eref{phinfarest} and \eref{phinnormest} follow from \eref{genlockernest} and a straightforward application of \eref{farptest} and \eref{totalnormest} respectively.
\end{Proof}

Corresponding to the formal kernel in \eref{formalkerndef},
we have the formal operator:
\be\label{formalsummopdef}
\sigma_n(\nu;\{\psi_j\}, \{\tilde{\psi}_j\}, H;f,x):=\int_\XX \Phi_n(\{\psi_j\}, \{\tilde{\psi}_j\}, H;x,y)f(y)d\nu(y),
\qquad f\in L^1,\ n>0,\ x\in\XX.
\ee
It is convenient to define
\be\label{sigmaenddef}
\sigma_0(\nu;\{\psi_j\}, \{\tilde{\psi}_j\}, H;f,x):=\int_\XX \Phi_0(\{\psi_j\}, \{\tilde{\psi}_j\},H;x,y)f(y)d\nu(y), \qquad f\in L^1,\  x\in\XX.
\ee
The following proposition lists some norm estimates for these operators.
\begin{prop}\label{gensigmaopprop}
We assume the set up as in Proposition~\ref{acttauberprop}. Let $\nu_1$ be a $d_1$-regular measure,  $\nu_2$ be a $d_1$-regular measure, and $1\le p\le r \le \infty$. Then for $f\in L^p(\nu_1)$
\bea\label{gensigmaopbd}
\|\sigma_n(\nu_1;\{\psi_j\}, \{\tilde{\psi}_j\}, H;f)\|_{\nu_2;r}&\le& 
c\left((c_1+(nd_1)^q)\tn\nu_1\tn_{R,d_1}\right)^{1/p'}\nonumber\\
&&\quad\quad
\times\left((c_1+(nd_2)^q)\tn\nu_2\tn_{R,d_2}\right)^{1/r} N^{Q-q+q(1/p-1/r)}\|f\|_{\nu_1;p}.\nonumber\\
\eea
\end{prop}

\begin{Proof}\ 
In this proof, we will abbreviate $\Phi_n(\{\psi_j\}, \{\tilde{\psi}_j\}, H; x,y)$ by $\Phi_n(x,y)$ and $\sigma_N(\nu_1;\{\psi_j\}, \{\tilde{\psi}_j\}, H;f)$ by $\sigma_n(f)$. Without loss of generality, we may assume also that $V_{S,Q}(H)=1$, and $\nu_j$ are positive measures. Using H\"older inequality and \eref{phinnormest}, we deduce that for $x\in\XX$,
\begin{eqnarray*}
|\sigma_n(f,x)|&\le& \int_\XX |\Phi_n(x,y)||f(y)|d\nu_1(y)\le \|\Phi_n(x,\circ)\|_{\nu_1;p'}\|f\|_{\nu_1;p}\\
& \le& c\left(c_1+(nd_1)^q)\tn\nu_1\tn_{R,d_1}\right)^{1/p'}n^{Q-q/p'}\|f\|_{\nu_1;p};
\end{eqnarray*}
i.e.,
\be\label{pf1eqn1}
\|\sigma_n(f)\|_{\nu_2;\infty} \le c\left((c_1+(nd_1)^q)\tn\nu_1\tn_{R,d_1}\right)^{1/p'}n^{Q-q/p'}\|f\|_{\nu_1;p}.
\ee
This proves \eref{gensigmaopbd} when $r=\infty$.

In particular, with $p=\infty$,
\be\label{pf1eqn2}
\|\sigma_n(f)\|_{\nu_2;\infty} \le c(c_1+(nd_1)^q)\tn\nu_1\tn_{R,d_1}n^{Q-q}\|f\|_{\nu_1;\infty}.
\ee
Switching the roles of $\{\psi_j\}$ and $\{\tilde{\psi}_j\}$ in \eref{phinnormest} (used with $\nu_2$ in place of $\nu$, $r=1$), the estimate becomes
$$
\int_\XX |\Phi_n(x,y)|d\nu_2(x)\le c(c_1+(nd_2)^q)\tn\nu_2\tn_{R,d_2}n^{Q-q}.
$$
Therefore,
\bea\label{pf1eqn3}
\|\sigma_n(f)\|_{\nu_2;1}&\le& \int_\XX\int_\XX|\Phi_n(x,y)||f(y)|d\nu_1(y)d\nu_2(x) = \int_\XX\left\{\int_\XX|\Phi_n(x,y)|d\nu_2(x)\right\}|f(y)|d\nu_1(y)\nonumber\\
& \le& c(c_1+(nd_2)^q)\tn\nu_2\tn_{R,d_2}n^{Q-q}\|f\|_{\nu_1;1}.
\eea
In view of the Riesz-Thorin interpolation theorem \cite[Theorem~1.1.1]{bergh}, \eref{pf1eqn2} and \eref{pf1eqn3} lead to
\be\label{pf1eqn4}
\|\sigma_n(f)\|_{\nu_2;p} \le c\left((c_1+(nd_1)^q)\tn\nu_1\tn_{R,d_1}\right)^{1/p'}
\left((c_1+(nd_2)^q)\tn\nu_2\tn_{R,d_2}\right)^{1/p}n^{Q-q}\|f\|_{\nu_1;p}, \qquad 1\le p\le\infty.
\ee
Hence, using \eref{pf1eqn1}, we obtain for $1\le r<r<\infty$,
\begin{eqnarray*}
\int_\XX |\sigma_n(f,x)|^rd\nu_2(x) &=&\int_\XX |\sigma_n(f,x)|^{r-p}|\sigma_n(f,x)|^pd\nu_2(x)\le \|\sigma_n(f)\|_{\nu_2;\infty}^{r-p}\|\sigma_n(f)\|_{\nu_2;p}^p\\
& \le&  c\left((c_1+(nd_1)^q)\tn\nu_1\tn_{R,d_1}\right)^{(r-p)/p'}
\left((c_1+(nd_1)^q)\tn\nu_1\tn_{R,d_1}\right)^{p/p'}\\
&& \qquad\qquad
\left((c_1+(nd_2)^q)\tn\nu_2\tn_{R,d_2}\right)
n^{(Q-q/p')(r-p)+(Q-q)p}\|f\|_{\nu_1;p}^r.
\end{eqnarray*}
This leads to \eref{gensigmaopbd} when $r<\infty$, and completes the proof.
\end{Proof}

In the sequel, we fix an infinitely differentiable, even function $h:\RR\to\RR$, nonincreasing on $[0,\infty)$, such that $h(t)=1$ if $|t|\le 1/2$, and $h(t)=0$ if $|t| \ge 1$.
We will write $g(t)=h(t)-h(2t)$. The following identities will be
used often without reference: For integers $n, k\ge 0$,
\be\label{hgidentities}
h\left(\frac{k}{2^n}\right)=h(k)+\sum_{j=1}^n g\left(\frac{k}{2^j}\right),
\qquad h\left(\frac{k}{2^n}\right)+\sum_{j=n+1}^\infty
g\left(\frac{k}{2^j}\right) =1.
\ee
We summarize the localization properties of three kernels which we will use in our proofs. Let $b$ be a mask of type $\beta\in\RR$. In the sequel, if $n>0$,
we will write $b_n(t)=b(nt)$.
\begin{prop}\label{kernprop}
Let $x, y\in\XX$, $n, N\ge 1$, $k=1,\cdots, R$, $d>0$, $\nu$ be a $d$-regular measure, $1\le p\le r\le \infty$. We have
\be\label{kernlocest}
| \Phi_n(h;x,y)|\le c\frac{n^{q}}{\max(1, (n\rho(x,y))^S)}, \quad \|\Phi_n(h;x,\circ)\|_{\nu;p}\le c\left((c_1+(nd)^q)\tn\nu\tn_{R,d}\right)^{1/p}n^{q/p'}.
\ee
\be\label{opkernlocest}
| L_{k,1}\Phi_n(h;x,y)|\le c\frac{n^{q+q_k}}{\max(1, (n\rho(x,y))^S)}, \quad \|L_{k,1}\Phi_n(h;x,\circ)\|_{\nu;p}\le c\left((c_1+(nd)^q)\tn\nu\tn_{R,d}\right)^{1/p}n^{q_k-q/p'}.
\ee
and
\bea\label{opgkernlocest}
|L_{k,1}\Phi_n(gb_N;x,y)| &\le& cN^{-\beta}\frac{n^{q+q_k}}{\max(1, (n\rho(x,y))^S)}, \nonumber\\
 \|L_{k,1}\Phi_n(gb_N;x,\circ)\|_{\nu;p}&\le& c\left((c_1+(nd)^q)\tn\nu\tn_{R,d}\right)^{1/p}n^{q_k-q/p'}N^{-\beta}.
\eea
Further, with the set up as in  Proposition~\ref{gensigmaopprop}, and writing $C_1=(c_1+(nd_1)^q)\tn\nu_1\tn_{R,d_1}$ and $
C_2=(c_1+(nd_2)^q)\tn\nu_2\tn_{R,d_2}$, we have
\be\label{sigmaopbd}
\|\sigma_n(\nu_1;h;f)\|_{\nu_2;r}\le  
cC_1^{1/p'}C_2^{1/r} n^{q(1/p-1/r)}\|f\|_{\nu_1;p}, 
\ee
\be\label{lksigmaopbd} 
\|L_k\sigma_n(\nu_1;h;f)\|_{\nu_2;r}\le  
cC_1^{1/p'}C_2^{1/r} n^{q_k+q(1/p-1/r)}\|f\|_{\nu_1;p},
\ee
and
\be\label{gbsigmaopbd}
\|L_k\sigma_n(\nu_1;gb_N;f)\|_{\nu_2;r}\le  
cC_1^{1/p'}C_2^{1/r} n^{q_k+q(1/p-1/r)}N^{-\beta}\|f\|_{\nu_1;p}.
\ee
\end{prop}

\begin{Proof}\ 
In view of Proposition~\ref{acttauberprop}, the first estimate in \eref{kernlocest} (respectively, \eref{opkernlocest})
follows  from \eref{vsqest} and \eref{gaussianbd} (respectively,  \eref{diffgaussianbd}). 
The second estimate in \eref{kernlocest}  follows from \eref{phinnormest} with the choices $\tilde{\psi}_j=\psi_j=\phi_j$, $H=h$, by observing from \eref{gaussianbd} that $Q=q$, and from \eref{vsqest} that $V_{S,q}(H)\le c$. Proposition~\ref{gensigmaopprop} yields  \eref{sigmaopbd} with the same choices.

The second estimate in \eref{opkernlocest} follows similarly, except with the choice $\tilde{\psi_j}=L_k\phi_j$ and the observation that \eref{diffgaussianbd} implies that $Q=q+q_k$. Proposition~\ref{gensigmaopprop} yields  \eref{lksigmaopbd} with the same choices.

It is easy to verify by induction that 
$$
\left|t^k\frac{d^k}{dt^k} ((1+t)^\beta
b(t))\right|=\left|t^k\frac{d^k}{dt^k}F_b(\log t)\right|\le c(b)c_2,
\qquad t> 0, \ k=0,\cdots,S,
$$
and hence, for $N\ge 1$,
\be\label{bderupbd}
\left|t^k\frac{d^k}{dt^k} ((1/N+t)^\beta b_N(t))\right|\le c(b)c_2N^{-\beta},
\qquad t> 0, \ k=0,\cdots,S+1.
\ee
Since $b(t)^{-1}$ is a mask of type $-\beta$, we record that
\be\label{bderlowbd}
\left|t^k\frac{d^k}{dt^k}((1/N+t)^\beta b_N(t))^{-1}\right|\le
c(b)c_2N^{\beta}, \qquad t> 0,\ k=0,\cdots,S+1.
\ee
Finally, if $g:\RR\to\RR$ is any compactly supported,
$S$ times continuously differentiable function, such that $g(t)=0$
on some neighborhood of $0$ then \eref{bderupbd}, \eref{bderlowbd} imply 
\be\label{gtimesbineq}
\||gb_N\||_S\le c(b,g)N^{-\beta}, \qquad \||g/b_N\||_S\le c(b,g)N^\beta,
\qquad N\ge 1.
\ee
In particular, \eref{vsqest}, and \eref{gtimesbineq} imply that $V_{S,Q}(gb_N)\le cN^{-\beta}$.
The first estimate in \eref{opgkernlocest} follows from Proposition~\ref{acttauberprop}. The second estimate follows similarly to the second estimate in \eref{opkernlocest}. The estimate \eref{gbsigmaopbd} follows from Proposition~\ref{gensigmaopprop} as with \eref{lksigmaopbd}.
\end{Proof}

\bhag{Integral representation}\label{intrepsect}
The objective of this section is to prove

\begin{theorem}\label{lkreptheo}
Let $1\le p\le \infty$, $\beta>\max_{1\le k\le m}q_k+q/p$,
and $f\in W_{p,b}$. Then for every $x\in \XX$,
\be\label{lkreprodkern}
L_k(f,x)=\int_\XX L_{k,1}G(x,y)\derf{f}{b}(y)d\mu^*(y).
\ee
Moreover, for $k=1,\cdots,R$,
\be\label{simulapprox}
\|f-\sigma_n(h;f)\|_\infty \le cn^{-q(\beta-1/p)}\|\derf{f}{b}\|_p,
\quad  \|L_kf-L_k\sigma_n(h;f)\|_\infty \le
cn^{-q(\beta-q_k-1/p)}\|\derf{f}{b}\|_p.
\ee
\end{theorem}

Our proof of Theorem~\ref{lkreptheo} requires  two inequalities: the Bernstein inequality and the Nikolskii inequality.

The Bernstein ineqality is the following.

\begin{lemma}\label{bernineqlemma}
Let  $1\le k\le R$. Then
\be\label{bernsteinineq}
\|L_k(P)\|_p \le c n^{q_k}\|P\|_p, \qquad P\in \Pi_n, \ n\ge 1.
\ee
\end{lemma}

\begin{Proof}\ 
If $P\in \Pi_n$, then a straightforward calculation using the orthogonality of $\{\phi_k\}$ and the facts
that $h(t)=1$ if $|t|\le 1/2$ and $\hat P(k)=0$ if $k>n$ shows that for $x\in\XX$,
\be\label{pf2eqn1}
P(x)=\int_\XX \Phi_{2n}(h;x,y)P(y)d\mu^*(y)=\sigma_{2n}(\mu^*;h,P), \qquad x\in \XX.
\ee
Therefore, $L_kP=L_k\sigma_{2n}(\mu^*;h;P)$. 
We now apply \eref{lksigmaopbd} with $2n$ in place of $n$,  $\nu_2=\nu_1=\mu^*$ (so that $d_1=d_2=0$, $\tn\nu_j\tn_{d_j}=1$), and $r=p$. 
This leads to \eref{bernsteinineq}.
\end{Proof}

The Nikolskii inequality is the following.

\begin{lemma}\label{nikollemma}
If $n>0$, $P\in\Pi_n$, $1\le p<r\le \infty$, then
\be\label{nikol}
\|P\|_r \le cn^{q(1/p-1/r)}\|P\|_p.
\ee
\end{lemma}
\begin{Proof}\ 
We  apply \eref{sigmaopbd} with $2n$ in place of $n$,  $\nu_2=\nu_1=\mu^*$ (so that $d_1=d_2=0$, $\tn\nu_j\tn_{d_j}=1$), and $P$ in place of $f$ to obtain
$$
\|\sigma_{2n}(\mu^*;h;P)\|_r\le cn^{q(1/p-1/r)}\|P\|_p.
$$
The estimate \eref{nikol} follows from \eref{pf2eqn1}.
\end{Proof}

We are now ready to prove Theorem~\ref{lkreptheo}.\\

\noindent\textsc{Proof of Theorem~\ref{lkreptheo}.}
In this proof only, let  $L$ denote any of the operators $L_k$,
and $A$ be the corresponding $q_k$ as defined in \eref{bernsteinineq}.
We observe using \eref{wbspacedef} that
$$
\sigma_{2^j}(gb_{2^j}; \derf{f}{b})=\sigma_{2^j}(g;f).
$$
In view of \eref{gtimesbineq}, we deduce that for integers $j\ge 1$,
\be\label{pf3eqn1}
\|\sigma_{2^j}(g;f)\|_p = \|\sigma_{2^j}(gb_{2^j}; \derf{f}{b})\|_p \le
c2^{-j\beta}\|\derf{f}{b}\|_p.
\ee
Since $\sigma_{2^j}(g;f)\in \Pi_{2^j}$, \eref{nikol} implies that
\be\label{pf3eqn2}
\|\sigma_{2^j}(g;f)\|_\infty \le c2^{-j(\beta-1/p)}\|\derf{f}{b}\|_p.
\ee
Hence, \eref{bernsteinineq} shows that
\be\label{pf3eqn3}
\|L\sigma_{2^j}(g; f)\|_\infty \le c2^{-j(\beta-A-1/p)}\|\derf{f}{b}\|_p.
\ee
Using the second identitiy in \eref{hgidentities}, we deduce
\eref{simulapprox} easily from \eref{pf3eqn2} and \eref{pf3eqn3}.
In particular, since $L$ is a closed operator, this shows that for $x\in \XX$,
$$
L(f,x)=\sum_{j=0}^\infty L\sigma_{2^j}(g;f,x) =\sum_{j=0}^\infty
L\sigma_{2^j}(gb_{2^j}; 
\derf{f}{b}, x),
$$
with the convergence being uniform. This leads to \eref{lkreprodkern}.
\qed

\bhag{An approximation result}\label{approxsect}

The purpose of this section is to prove the following result, analogous to \cite[Theorem~6.3]{bdint}.

\begin{theorem}\label{networktheo}
Let $1\le p\le \infty$, $\beta>\max_{1\le k \le R}q_k +q/p'$, and
\be\label{networkdef}
\Psi(x)=\sum_{k=1}^R\sum_{j=1}^M a_{k,j}L_{k,1}G(y_j,x), \qquad x\in\XX.
\ee
There exists an integer $N^*\sim\eta^{-1}$ such that for $n\ge N^*$,
\be\label{networkapprox}
\|\Psi-\sigma_n(h;\Psi)\|_p \le \frac{1}{2}\|\Psi\|_p.
\ee
\end{theorem}

As in \cite{bdint}, the main problem is to relate $\|\Psi\|_p$ and
the coefficients $c_{k,j}$, via $\sigma_{2^m}(g;\Psi)$.
The main technical problem is the following.
The mapping $y\mapsto \sigma_{2^m}(g;L_{k,1}G(y,x))$ is not in $\Pi_{2^m}$.
So, we cannot use  an analogue of \cite[Proposition~6.3]{bdint}
to obtain a full estimate of the form \cite[Theorem~6.2(b)]{bdint}.
The following Proposition~\ref{coeffprop} (which we will call the coefficient inequalities) serves as a substitute.

In order to state the coefficient inequalities, we will use the following notations: ${\bf a}_k=(a_{k,1},\cdots,a_{k,j})$, and
\be\label{partialpsi}
\Psi_k(x)=\sum_{j=1}^M a_{k,j}L_{k,1}G(y_j,x), \qquad x\in\XX,
\ee
so that $\Psi=\sum_{k=1}^R\Psi_k$. For any sequence ${\bf d}$,
$$
\|{\bf d}\|_p :=\left\{\begin{array}{ll}
\left\{\sum_{j=1}^\infty |d_j|^p\right\}^{1/p}, &\mbox{ if $1\le p<\infty$},\\
\sup_{1\le j\le \infty}|d_j|, &\mbox{ if $p=\infty$,}
\end{array}\right.
$$
with a Euclidean vector extended to a sequence by padding with $0$'s. The coefficient inequalities are given in the following proposition.

\begin{prop}\label{coeffprop}
Let $1\le p\le \infty$, and $1\le k \le R$.\\
{\rm (a)} For integer $m$ with $2^m\eta\ge 1$, we have
\be\label{partpsibycoeff}
\|\sigma_{2^m}(g;\Psi_k)\|_p \le c2^{m(q_k-\beta+q/p')}\|{\bf a}_k\|_p.
\ee
{\rm (b)} We have
\be\label{coeffbyfullnet}
\|{\bf a}_k\|_p \le c\eta^{q_k-\beta+q/p'}\|\Psi\|_p.
\ee
{\rm (c)} For integer $m$ with $2^m\eta\ge 1$, we have
\be\label{fullpsicoeff}
\|\sigma_{2^m}(g;\Psi)\|_p \le c\sum_{k=1}^R
(2^m\eta)^{q_k-\beta+q/p'}\|\Psi\|_p,
\ee
and
\be\label{fullpsinormproperty}
\sum_{k=1}^R \eta^{-q_k}\|{\bf a}_k\|_p \le c\eta^{q/p'-\beta}\|\Psi\|_p.
\ee
\end{prop}

\begin{Proof}\ 
The proof of part (a) mimics that of \cite[Theorem~6.2(a)]{bdint}.
We observe that for $x\in\XX$,
$$
\sigma_{2^m}(g;L_{k,1}G(y_j,\circ), x) =
\sum_{i}
g(\lambda_i/2^m)b(\lambda_i)L_k(\phi_i)(y_j)\phi_i(x)=L_{k,1}\Phi_{2^m}(gb_{2^m};y_j,x).
$$
Using \eref{opgkernlocest} with  $p=1$, $\nu=\mu^*$ (so that $d=0$, $\tn\nu\tn_{R,d}=1$),  we obtain
\be\label{pf4eqn1}
\|\sigma_{2^m}(g;L_{k,1}G(y_j,\circ)\|_1\le c2^{m(q_k-\beta)}.
\ee
Therefore,
\be\label{pf4eqn2}
\|\sigma_{2^m}(g; \Psi_k)\|_1\le \sum_{j=1}^M
|a_{k,j}|\|\sigma_{2^m}(g;L_{k,1}G(y_j,\circ))\|_1\le
c2^{m(q_k-\beta)}\|{\bf a}_k\|_1.
\ee

Next, we use \eref{opgkernlocest} again with $2^m\ge \eta^{-1}$ in place of $n$, $p=1$ and $\nu$ to be the measure $\tau$ as defined in Lemma~\ref{minsepmeaslemma} (so that $d=\eta$, $\tn\tau\tn_{R,\eta}\le c$) to deduce that
\bea\label{pf4eqn3}
\left\|\sum_{j=1}^M |\sigma_{2^m}(g;L_{k,1}G(y_j,\circ))|\right\|_\infty
&=&\left|\eta^{-q}\int_\XX |\sigma_{2^m}(g;L_{k,1}G(\circ,y))|d\tau(y)\right\|_\infty\nonumber\\
& \le& c2^{m(q+q_k-\beta)}(1+(2^m\eta)^q)(2^m\eta)^{-q} \le
c2^{m(q+q_k-\beta)}.
\eea
Therefore,
\bea\label{pf4eqn5}
\|\sigma_{2^m}(g; \Psi_k)\|_\infty&\le& \sum_{j=1}^M
|a_{k,j}|\|\sigma_{2^m}(g;L_{k,1}G(y_j,\circ)\|_\infty\nonumber\\
&\le& \|{\bf a}_k\|_\infty\left\|\sum_{j=1}^M
|\sigma_{2^m}(g;L_{k,1}G(y_j,\circ))|\right\|_\infty \le
c2^{m(q+q_k-\beta)}\|{\bf a}_k\|_\infty.
\eea
The estimate \eref{partpsibycoeff} follows from \eref{pf4eqn2}
and \eref{pf4eqn5} by an application of the
Riesz--Thorin interpolation theorem \cite[Theorem~1.1.1]{bergh}
to the operator ${\bf a}_k\mapsto \Psi_k$. 

To prove part (b), we fix $k$, and find ${\bf d}\in\RR^M$ such that 
\be\label{pf4eqn8}
\langle {\bf a}_k, {\bf d}\rangle = \|{\bf a}_k\|_p, \qquad \|{\bf d}\|_{p'}=1.
\ee
We then recall the functions $\phi_{k,j}$ defined in the assumptions on
$L_k$, and set 
$$
F(x)=\sum_{j=1}^M d_j\phi_{k,j}^{(b)}(x), \qquad x\in \XX.
$$
We estimate $\|F\|_{p'}$. We consider the case when $1\le p' <\infty$,
the case when $p'=\infty$ is only simpler.
Since each $\phi_{k,j}^{(b)}$ is supported on a neighborhood of $y_j$
with diameter $\le \eta/3$, for any given $x\in\XX$,
there is at most one $j'$ such that $\phi_{k,j'}^{(b)}(x)\not=0$. Then 
$$
|F(x)|^{p'}=|d_{j'}|^{p'}|\phi_{k,j'}^{(b)}(x)|^{p'},
$$
and thus, for any $x\in\XX$,
\be\label{pf4eqn6}
|F(x)|^{p'} = \sum_{j=1}^M |d_j|^{p'}|\phi_{k,j}^{(b)}(x)|^{p'}.
\ee
Since each $\phi_{k,j}^{(b)}$ is supported on a neighborhood of $y_j$
with diameter $\le \eta/3$, \eref{dualfnnormest} and
\eref{ballmeasurecond} shows that
$$
\int_\XX |\phi_{k,j}^{(b)}(x)|^{p'}d\mu^*(x) \le \int_{\BB(y_j,\eta/3)}
|\phi_{k,j}^{(b)}(x)|^{p'}d\mu^*(x) \le c\eta^{(q_k-\beta)p'+q}.
$$
Consequently, \eref{pf4eqn6} implies that
\be\label{pf4eqn7}
\|F\|_{p'} \le c\eta^{q_k-\beta+q/p'}.
\ee
Next, using Theorem~\ref{lkreptheo} and \eref{dualfndef} from
the assumptions on $L_k$, we obtain that
for any $k'=1,\cdots,R$, 
\begin{eqnarray*}
\lefteqn{\int_\XX L_{k',1}G(y_{j'},x)F(x)d\mu^*(x) =\sum_{j=1}^M d_j
\int_\XX L_{k',1}G(y_{j'},x)\phi_{k,j}^{(b)}(x)d\mu^*(x)}\\
&=& \sum_{j=1}^M d_j L_{k'}(\phi_{k,j})(y_{j'}) =\left\{\begin{array}{ll}
d_{j'}, &\mbox{ if $k'=k$,}\\
0, &\mbox{ otherwise.}
\end{array}\right.
\end{eqnarray*}
Consequently, the definition \eref{pf4eqn8} of ${\bf d}$ shows that
$$
\int_\XX \Psi(x)F(x)d\mu^*(x) = \sum_{k'=1}^R \sum_{j'=1}^M
a_{k',j'}\int_\XX
L_{k',1}G(y_{j'},x)F(x)d\mu^*(x)=\sum_{j'=1}^M a_{k,j'}d_{j'} = \|{\bf a}_k\|_p.
$$
Together with \eref{pf4eqn7}, this shows that
$$
\|{\bf a}_k\|_p =\int_\XX \Psi(x)F(x)d\mu^*(x) \le \|\Psi\|_p\|F\|_{p'}
\le c\eta^{q_k-\beta+q/p'}\|\Psi\|_p.
$$
This proves part (b).

Part (c) is a straightforward consequence of parts (a) and (b)
and the fact that $\Psi=\sum_{k=1}^R \Psi_k$.
\end{Proof}

\noindent
\textsc{Proof of Theorem~\ref{networktheo}.} In this proof,
let $A=\max_{1\le k\le R}q_k$.
In light of \eref{hgidentities} and \eref{fullpsicoeff}
we obtain for $N\ge \eta^{-1}$, $\beta> A+q/p'$ that
$$
\|\Psi-\sigma_{2^N}(h;\Psi)\|_p \le
\sum_{m=N+1}\|\sigma_{2^m}(g;\Psi)\|_p \le c\sum_{k=1}^R\sum_{m=N+1}
(2^m\eta)^{q_k-\beta+q/p'}\|\Psi\|_p \le
c_1(2^N\eta)^{A-\beta+q/p'}\|\Psi\|_p.
$$
If $n\ge 2^{N+1}$, then
$$
\|\Psi-\sigma_n(h;\Psi)\|_p \le c\inf_{P\in\Pi_{n/2}}\|\Psi-P\|_p
\le c\inf_{P\in\Pi_{2^N}}\|\Psi-P\|_p\le
c\|\Psi-\sigma_{2^N}(h;\Psi)\|_p\le c_2(2^N\eta)^{A-\beta+q/p'}\|\Psi\|_p.
$$
We may choose $N$ so that $2^N\sim \eta^{-1}$ and the rightmost
expression above is $\le (1/2)\|\Psi\|_p$.
Then we have proved \eref{networkapprox} with $N^*=2^{N+1}$. \qed
 
\bhag{Proofs of the main results}\label{pfsect}

\textsc{Proof of Theorem~\ref{golombtheo}.}\\

We define the inner product
$$
\ip{g_1}{g_2}=\sum_{k=0}^\infty \frac{\hat{g_1}(k)\hat{g_2}(k)}{b(\lambda_k)^2}.
$$
We wish to  show first that there exist the coefficients $a_{k,j}$ such that
\be\label{kerninterp}
L_iP(y_\ell)=f_{i,\ell}, \qquad \ell=1,\cdots,M_n, \ i=1,\cdots,R,
\ee
and with this choice, $P$ is the solution of the minimization
problem \eref{golomb}.

We observe that 
$$
\sum_{k,j}\sum_{i,\ell}d_{k,j}d_{i,\ell}L_{i,2}L_{k,1}\GG(y_j,y_\ell)
=\sum_{\nu=0}^\infty
b(\lambda_\nu)^2\left\{\sum_{k,j}d_{k,j}L_k(\phi_\nu)(y_j)\right\}^2 \ge 0.
$$
If the infinite sum is equal to $0$, then
$$
\sum_{k,j}d_{k,j}L_k(\phi_\nu)(y_j) =0,\qquad \nu=0,1,\cdots.
$$
Consequently, for every $x\in\XX$,
$$
0=\sum_{\nu=0}^\infty
b(\lambda_\nu)^2\sum_{k,j}d_{k,j}L_k(\phi_\nu)(y_j)\phi_\nu(x)
=\sum_{k,j} d_{k,j}L_{k,1}\GG(y_j,x).
$$

Next, we observe by a straightforward calculation that for any $f$
in the domain of the $L_k$'s,
\be\label{ipopval}
\ip{L_{k,1}\GG(y,\circ)}{f} =L_k(f)(y).
\ee
So, for each $i$ and $\ell$,
$$
d_{i,\ell}=L_i(\phi_{i,\ell})(y_\ell)=\left\langle\sum_{k,j}
d_{k,j}L_{k,1}\GG(y_j,\circ), \phi_{i,\ell}\right\rangle=0.
$$ 
 Thus, the matrix of the system of equations \eref{kerninterp}
is positive definite, and hence, \eref{kernsol} has a unique
solution given by $P$.
 
 Let $g$ be another candidate for the minimization problem. Then
\begin{eqnarray*}
 \ip{P}{g-P}&=&\sum_{k,j}a_{k,j}\ip{L_{k,1}\GG(y_j,\circ)}{g-P} =
\sum_{k,j}a_{k,j}(L_k(g)(y_j)-L_k(P)(y_j))\\
 & =&  \sum_{k,j}a_{k,j}(L_k(f)(y_j)-L_k(f)(y_j))=0.
 \end{eqnarray*}
 Hence,
 $$
 \|\derf{g}{b}\|_2^2=\ip{g}{g} =\ip{g-P}{g-P}+\ip{P}{P}\ge \|\derf{P}{b}\|_2^2.
 $$
 This proves that $P$ is the solution of the minimization problem. 
\qed

\par
\noindent
\textsc{Proof of Theorem~\ref{feasibletheo}.}\\

This proof is almost verbatim the same as the proof of
\cite[Theorem~5.1]{bdint}. We will point out only the differences.
We will omit the notation $n$ in our proof as in the other paper.
Let ${\bf a}\in \RR^{RM}$ be a row--major ordering of the matrix
$(a_{k,j})_{k=1,\cdots,R, \ j=1,\cdots, M}$. We define
\be\label{pf5eqn1}
\tn{\bf a}\tn_{RM}^* := \left\| \sum_{k=1}^R\sum_{j=1}^M
a_{k,j}L_{k,1}G(y_j,\circ)\right\|_{p'}.
\ee
In view of \eref{fullpsinormproperty} this is a norm on $\RR^{RM}$.
Let $N^*$ be chosen so that Theorem~\ref{networktheo} holds with
$p'$ in place of $p$, and $D$ be the dimension of $\Pi_{N^*}$.
For ${\bf d}\in \RR^D$, we define
\be\label{pf5eqn2}
\tn {\bf d}\tn_D := \left\|\sum_{\lambda_j< N^*} b(\lambda_j)^{-1}
d_j\phi_j\right\|_p.
\ee
In place of $F^{(-s)}$ in the proof in \cite{bdint}, we need $F^{(1/b)}$.
In place of the matrix ${\bf A}$ in \cite{bdint}, we take the matrix
appropriate for the interpolation problem which we are dealing with;
i.e., a matrix indexed by $(k,j)$ and $m$ so that the
$((k,j),m)$--th entry is $L_k(\phi_m)(y_j)$.
The rest of the proof is the same as in \cite{bdint}
with obvious minor changes; e.g., $\XX$ in place of $[-\pi,\pi]^q$,
bivariate kernels in place of convolutions, etc.
\qed

We will prove Theorem~\ref{feasimpconvtheo} in much greater
generality for future reference in the form of
Theorem~\ref{convergetheo} below.

Let $\XX$ be a separable Banach space with norm $\|\circ\|$,
$\XX^*$ be its dual space with dual norm $\|\circ\|^*$.
Let $R\ge 1$ be an integer, and for each integer $n\ge 1$, $1\le k\le R$,
$x_{k,n}^*\in\XX^*$,  with $\|x_{k,n}^*\|^*\le 1$.
We assume that for each $k$, $x_{k,n}^*\rightarrow x_k^*$
in the weak-* topology.
Necessarily, each $x_k^*\in\XX^*$ and $\|x_k^*\|^* \le 1$. 

Let $\YY$ be another normed linear space with norm $|\circ|_\YY$,
continuously embedded, and hence, identified with a subspace of $\XX$.
We assume that the unit ball
$$
\mathcal{B} =\{g\in\YY : |g|_\YY \le 1\}
$$
is compact in $\XX$. Let $V_1\subseteq \cdots V_n\subseteq V_{n+1}
\subseteq \cdots$ be a sequence of subsets of $\YY$.

\begin{theorem}\label{convergetheo}
Let $f\in\YY$, and we assume that there exists $v_n\in V_n$ such that
\be\label{feasible}
|v_n|_\YY =\min\{|g|_\YY : g\in V_n,\ x_{k,n}^*(g)=x_k^*(f),
\ k=1,\cdots,R\}\le c|f|_\YY.
\ee
Then $x_k^*(v_n)\rightarrow x_k^*(f)$ as $n\to\infty$.
\end{theorem}

This theorem easily follows from the following lemma.
If $\delta>0$, and $K\subset\XX$ is a compact set, we say that a
subset $A$ of $K$ is $\delta$ separated if 
$$
\min_{f_1,f_2\in A\atop f_1\not= f_2}\|f_1-f_2\| \ge \delta.
$$
Since $K$ is compact, such a set is necessarily finite.

\begin{lemma}\label{weakstarlemma}
Let $\XX$ be a separable Banach space,  $\{y_m^*\}_{m=0}^\infty$
be a sequence in the unit ball of $\XX^*$ converging to $y^*\in\XX^*$
in the weak-* topology. If $K\subset \XX$ is a compact set, then
$$
\lim_{m\to\infty}\sup_{f\in K}|y_m^*(f)-y^*(f)|=0.
$$
\end{lemma}
\begin{Proof}\ 
Let $\e>0$. We consider the set $A$ to be the maximal $\e/3$
separated subset of $K$. If 
$$
f_0\in K\setminus\bigcup_{g\in A}\{f\in \XX : \|g-f\| \le \e/3\}.
$$
then we may add $f_0$ to $A$ to obtain a larger $\e/3$ separated
subset of $K$. 
Therefore,
\be\label{pf6eqn1}
K\subseteq \bigcup_{g\in A}\{f\in \XX : \|g-f\| \le \e/3\}.
\ee
In view of the weak-* convergence, we obtain an integer $N\ge 1$ such that
$$
\sup_{m\ge N,\ g\in A}|y_m^*(g)-y^*(g)|\le \e/3.
$$
Let $f\in K$. In view of \eref{pf1eqn1}, we obtain $g\in A$ such that
$$ 
\|f-g\|\le \e/3.
 $$
Using the fact that $\|y_m^*\|^*\le 1$, $\|y^*\|^*\le 1$,
we obtain for $m\ge N$,
$$
 |y_m^*(f)-y^*(f)|\le |y_m^*(f)-y_m^*(g)| +|y_m^*(g)-y^*(g)|+|y^*(g)-y^*(f)|
 \le \|f-g\| +\e/3 +\|f-g\| \le \e.
 $$
 This completes the proof.
\end{Proof}

\par
\noindent
\textsc{Theorem~\ref{feasimpconvtheo}.}
We apply Theorem~\ref{convergetheo} with the following choices.
The space $\XX$ is defined as follows.
We consider the space of all functions  $f\in L^p$ for which
$L_k(f)$ is well defined, with the norm
$$
\|f\|=\sum_{k=1}^R\|L_k(f)\|_p.
$$
Then $\XX$ is the closure of the space of all diffusion polynomials
in the sense of this norm.
The functionals are defined by
$$
x_{k,n}^*(f)=d_kL_k(f)(x_n), \quad x_k^*(f)=c_k L_k(f)(x_0),
$$
where the constants are chosen to bring the linear functionals
into the unit ball of $\XX^*$.
These will depend only on $L_k$ and $\XX$, not on the
individual points $x_n$. \qed

\bhag{Numerical Simulations}\label{num_results}
In this section we present numerical simulations of Birkhoff
interpolation. By noting the one-to-one correspondance
between even trigonometric polynomials and algebraic
polynomials of the same degree, we focus on 
interpolation on $[-1,1]$ and $[-1,1]^{2}$.
To show that this method works for a general basis,
we also include examples of interpolation on the sphere $\SS^{2}$
using real spherical harmonics.
Examples of both 1D and 2D interpolation 
are included,
but we emphasize 2D interpolation due to its challenging nature.
In two dimensions, we restrict ourselves to interpolating on subsets
of $[-1,1]^{2}$ and $\SS^{2}$.
In particular, we look at interpolating the function

\begin{equation}
f_{R}(x,y) = \frac{1}{1+ R(x^2 + y - 0.3)^{2}}
       + \frac{1}{1+ R(x + y - 0.4)^{2}}
       + \frac{1}{1+ R(x + y^{2} - 0.5)^{2}}
       + \frac{1}{1+ R(x^{2} + y^{2} - 0.25)^{2}},
\label{eq:runge_func}
\end{equation}

\noindent
with Runge functions on two parabolas, one line, and one circle
on $[-1,1]^{2}$.
Here, the $R$ parameter controls the difficulty of the interpolation,
with larger $R$ values corresponding to more difficult interplation problems.
Most of the tests use $R = 25$, although some examples
(Tables~\ref{tab:2d_error_tensor_float32_easy} and
\ref{tab:2d_error_tensor_easy})
use $R = 9$ in order to show that we approach machine precision even
in 2D interpolation.
In one dimension, we interpolate $f_{25}(x,-0.96)$.
On the sphere, we look at interpolating the function

\begin{equation}
g(x,y,z) = \frac{1}{1 + \left(\cos7x + \cos7y + \cos7z\right)^{2}},
\label{eq:runge_func_sphere}
\end{equation}

\noindent where we restrict $g$ to $\SS^{2}$.
Unless stated otherwise, the examples are computed using double precision
and the derivative information consists of directional derivatives
along the coordinate axes.
There does not appear to be any standard software packages for Birkhoff
interpolation when the point distribution is arbitrary, so we
are not able to compare our method with others.

For $p=2$, MSN Birkhoff interpolation reduces to solving
\begin{equation}
\min_{Va=f}||D_{s}a||_{2},
\end{equation}
where $D_{s}$ is a diagonal positive-definite matrix with
condition number $O(n^{s})$, $V$ is a Chebyshev-Vander\-monde
matrix, $a$ is the vector containing the Chebyshev interpolation
coefficients, and $f$ is the vector containing the function
and derivative values.
This linear system is solved in a similar way described in
\cite{bdint}.

In the one dimensional case, we interpolate $f_{25}(x,-0.96)$ on $[-1,1]$
using both single precision and double precision
on equally-spaced points. Results for function and derivative
information at equally spaced points are presented in
Tables~\ref{tab:1d_error_tensor_float32} and \ref{tab:1d_error_tensor}.
We obtain similar results in
Tables~\ref{tab:1d_error_tensor_interlace_float32}
and \ref{tab:1d_error_tensor_interlace}
when we interlace the function and derivative information.
In the 2D interpolation below,
the point $(-0.96,-0.96)$ was the point with the largest
error on the $21\times 21$ tensor grid for $s = 8$, and this
is why we chose to the 1D function to be $f_{25}(x,-0.96)$.
In both single and double precision
the interpolation error approaches machine epsilon with an appropriately
chosen $s$-value and increasing points. The error is computed by taking
the difference between the function and the approximation on
$10n$ equally spaced points, normalized by the maximum function
value.
The mesh norm is 

\begin{equation}
m = \left\lceil \frac{2\pi}{\min_{i\ne j}
\left|\cos^{-1}(x_{i}) - \cos^{-1}(x_{j})\right|
}
\right\rceil.
\end{equation}

\noindent and the interpolation order is taken to be $m$.

We showcase the true power of MSN Birkhoff interpolation by
using a variety of point distributions in two dimensions.
In particular, the interpolation of $f_{25}(x,y)$
is performed on a tensor
grid (Tables~\ref{tab:2d_error_tensor_float32}
and \ref{tab:2d_error_tensor}), on a
tensor grid intersected with an annulus
(Table~\ref{tab:2d_error_tensor_annulus}), and on a
annular grid (Table~\ref{tab:2d_error_annular_grid}).
The maximum error is computed by taking the maximum difference on a
$10n\times 10n$ grid for an $n\times n$ grid interpolation.
When we interpolate with function and derivative values 
interlaced, we see that also see convergence,
whether the directional derivatives are parallel to the coordinate
axes
(Tables~\ref{tab:2d_error_tensor_interlace_parallel_float32}
and \ref{tab:2d_error_tensor_interlace_parallel})
or not
(Tables~\ref{tab:2d_error_tensor_interlace_not_float32}
and \ref{tab:2d_error_tensor_interlace_not}).
Furthermore, MSN Birkhoff interpolation can use a variety
of bases.
Tables~\ref{tab:2d_error_sphere_float32} and
Table~\ref{tab:2d_error_sphere}
show results for interpolation on $\SS^{2}$
using real spherical harmonics.
Additionally, we present results in
Table~\ref{tab:2d_error_sphere_north_south}
for when points are clustered near
the north and south poles, requiring interpolation points
to have polar angle
$\theta\in[0,\frac{\pi}{3}]\cup[\frac{2\pi}{3},\pi]$.
We see a general trend of decreasing error in these regions.
By interpolating $f_{9}(x,y)$ on tensor grids on $[-1,1]^{2}$
(Tables~\ref{tab:2d_error_tensor_float32_easy} and
\ref{tab:2d_error_tensor_easy}),
we see that in 2D problems we still obtain errors that
approach machine precision.
The interpolation order for two-dimensional problems
is computed similar to the one-dimensional
case.

Larger $s$ values give greater derivative control and generally
more accurate results; however, large $s$ values
lead to high condition numbers, so care
must be taken to arrive at an accurate solution.
We use a variant of the  complete orthogonal decomposition of
\cite{Hough95completeorthogonal} because the ill-conditioning in our method
results primarily from a diagonal
matrix; also see \cite{bdint}
and references therein for more details of the numerical technique.
The difficulty in obtaining low interpolation error is
due to the long time required to run the dense matrix
computations. Future work will be devoted to investigating
and implementing fast algorithms arising from the structured equations.

We have not shown this method to be numerically stable; however,
all of the numerical examples are solved using the same algorithm
in both single and double precision. The one-dimensional data
indicate that we approach machine epsilon in both single and double
precision, while the two-dimensional data show that increasing
data points generally leads to decreased error.

\begin{figure}
\centering
\begin{subfigure}{0.45\textwidth}
\includegraphics[width=\textwidth]{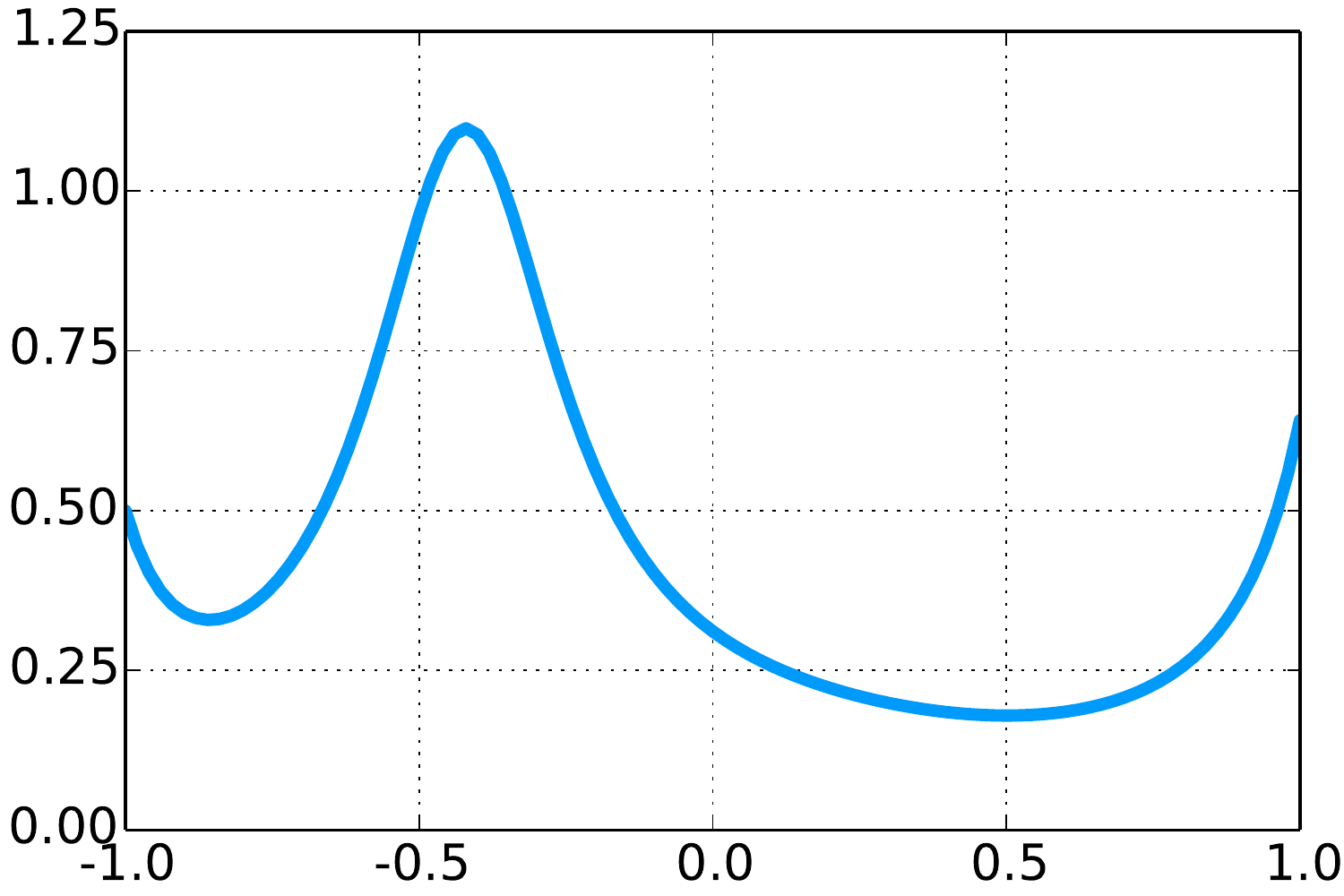}
\end{subfigure}
\begin{subfigure}{0.45\textwidth}
\includegraphics[width=\textwidth]{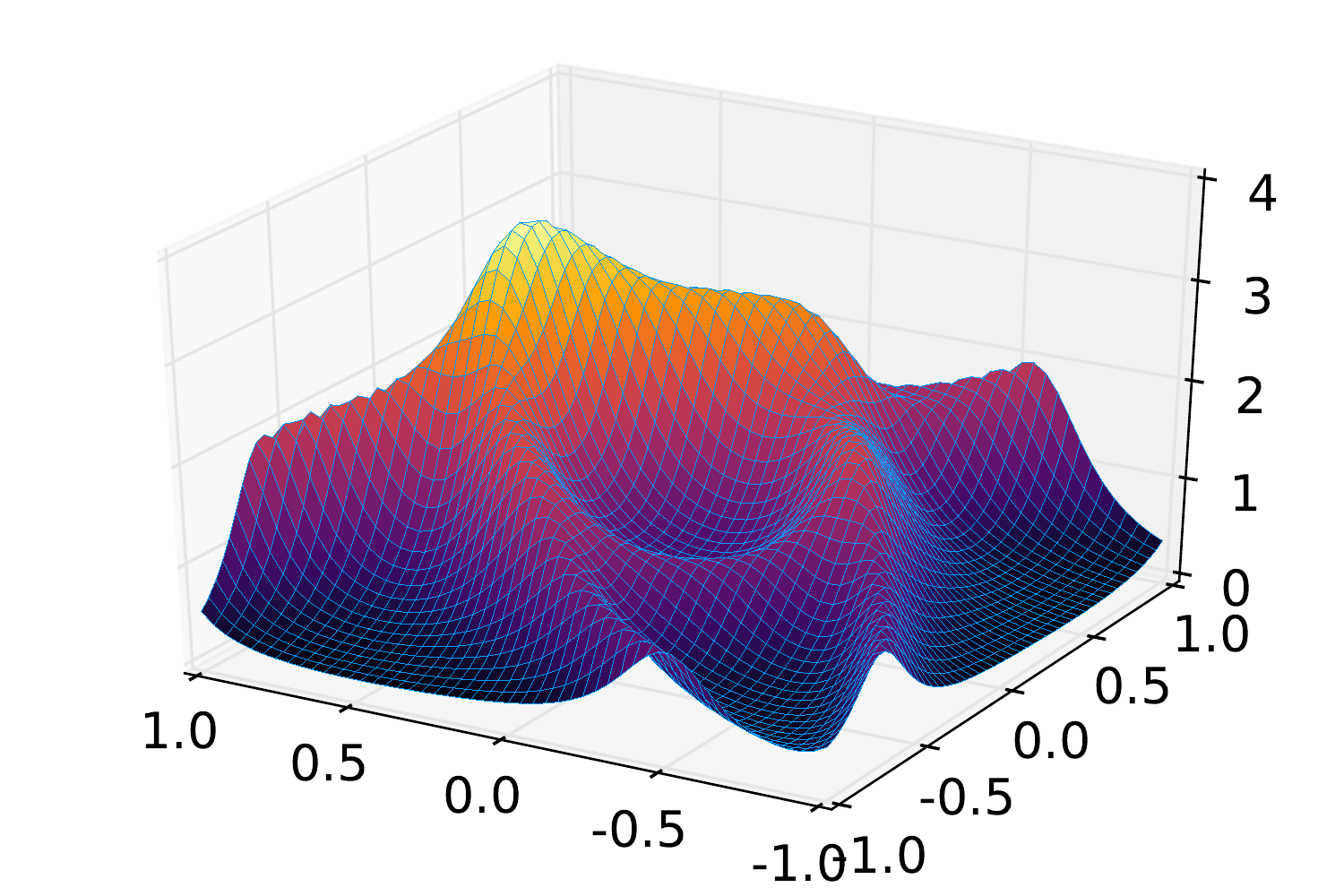}
\end{subfigure}
\caption{Test functions $f_{25}(x,-0.96)$ and $f_{25}(x,y)$ from
Eq.~(\ref{eq:runge_func}).
}
\end{figure}

\begin{table}
\centering
\begin{tabular}{|c||c|c|c|c|}
\hline
$n$ & $s = 2$ & $s = 3$ & $s = 4$ & $s = 5$ \\
\hline
\hline
11 & 1.40e-02 & 9.98e-03 & 1.93e-02 & 4.28e-02 \\
21 & 1.12e-03 & 9.15e-04 & 2.77e-04 & 1.01e-04 \\
31 & 1.87e-03 & 1.23e-04 & 3.07e-05 & 6.54e-06 \\
41 & 1.73e-03 & 4.77e-05 & 1.95e-06 & 1.66e-06 \\
51 & 1.47e-03 & 5.14e-05 & 4.29e-06 & 3.34e-06 \\
61 & 1.22e-03 & 4.30e-05 & 2.77e-06 & 2.99e-06 \\
\hline
\end{tabular}
\caption{
Interpolation error of MSN 1D Birkhoff interpolation
for various $s$ values on $[-1,1]$ for the function $f_{25}(x,-0.96)$.
Function and derivative information is given at $n$ equally spaced points.
Error has been normalized by
maximum function value.
All of the computations in this table were performed
in \textbf{single precision}.
}
\label{tab:1d_error_tensor_float32}
\end{table}

\begin{table}
\centering
\begin{tabular}{|c||c|c|c|c|c|c|}
\hline
$n$ & $s = 2$ & $s = 4$ & $s = 6$ & $s = 8$ & $s = 10$ & $s = 12$ \\
\hline
\hline
21 & 1.12e-03 & 2.78e-04 & 1.77e-04 & 2.57e-03 & 1.83e-02 & 7.56e-02 \\
41 & 1.73e-03 & 1.86e-06 & 5.32e-07 & 5.70e-08 & 1.23e-07 & 3.18e-07 \\
61 & 1.22e-03 & 2.73e-06 & 1.65e-07 & 1.49e-08 & 1.95e-09 & 5.15e-10 \\
81 & 8.37e-04 & 1.31e-06 & 3.84e-08 & 2.11e-09 & 1.51e-10 & 2.70e-11 \\
101 & 6.04e-04 & 6.30e-07 & 1.05e-08 & 3.68e-10 & 1.36e-11 & 3.76e-12 \\
121 & 4.55e-04 & 3.27e-07 & 3.35e-09 & 7.92e-11 & 1.42e-12 & 4.42e-12 \\
141 & 3.54e-04 & 1.83e-07 & 1.21e-09 & 2.00e-11 & 3.71e-13 & 4.32e-12 \\
161 & 2.83e-04 & 1.09e-07 & 4.78e-10 & 5.87e-12 & 7.24e-13 & 1.62e-11 \\
\hline
\end{tabular}
\caption{
Interpolation error of MSN 1D Birkhoff interpolation
for various $s$ values on $[-1,1]$ for the function $f_{25}(x,-0.96)$.
Function and derivative information is given at $n$ equally spaced points.
Error has been normalized by
maximum function value.
}
\label{tab:1d_error_tensor}
\end{table}

\begin{table}
\centering
\begin{tabular}{|c||c|c|c|c|}
\hline
$n$ & $s = 2$ & $s = 3$ & $s = 4$ & $s = 5$ \\
\hline
\hline
11 & 1.48e-02 & 1.25e-02 & 1.69e-02 & 2.50e-02 \\
21 & 2.07e-03 & 1.08e-03 & 6.39e-04 & 5.74e-04 \\
31 & 1.68e-03 & 1.59e-04 & 6.03e-05 & 3.92e-05 \\
41 & 1.36e-03 & 4.69e-05 & 7.74e-06 & 4.02e-06 \\
51 & 1.10e-03 & 3.13e-05 & 2.86e-06 & 3.03e-06 \\
61 & 8.93e-04 & 2.12e-05 & 4.02e-06 & 3.80e-06 \\
\hline
\end{tabular}
\caption{
Interpolation error of MSN 1D Birkhoff interpolation
for various $s$ values on $[-1,1]$ for the function $f_{25}(x,-0.96)$.
Function information is given at $n$ equally spaced points, while
derivative information are given at $n-1$ points in between.
Error has been normalized by
maximum function value.
All of the computations in this table were performed
in \textbf{single precision}.
}
\label{tab:1d_error_tensor_interlace_float32}
\end{table}

\begin{table}
\centering
\begin{tabular}{|c||c|c|c|c|c|c|}
\hline
$n$ & $s = 2$ & $s = 4$ & $s = 6$ & $s = 8$ & $s = 10$ & $s = 12$ \\
\hline
\hline
21 & 2.07e-03 & 6.41e-04 & 6.34e-04 & 1.25e-03 & 8.21e-03 & 5.86e-02 \\
41 & 1.36e-03 & 6.89e-06 & 1.94e-06 & 1.33e-06 & 1.06e-06 & 9.37e-07 \\
61 & 8.93e-04 & 1.30e-06 & 4.58e-08 & 7.60e-09 & 2.47e-09 & 1.88e-09 \\
81 & 6.07e-04 & 3.98e-07 & 1.86e-09 & 4.37e-11 & 8.73e-11 & 4.54e-11 \\
101 & 4.31e-04 & 1.55e-07 & 1.21e-09 & 7.21e-11 & 1.55e-11 & 3.27e-12 \\
121 & 3.22e-04 & 7.13e-08 & 5.72e-10 & 2.89e-11 & 3.25e-12 & 5.29e-12 \\
141 & 2.48e-04 & 3.68e-08 & 3.12e-10 & 1.07e-11 & 7.22e-13 & 7.30e-13 \\
161 & 1.97e-04 & 2.11e-08 & 1.68e-10 & 4.17e-12 & 1.70e-13 & 1.95e-12 \\
\hline
\end{tabular}
\caption{
Interpolation error of MSN 1D Birkhoff interpolation
for various $s$ values on $[-1,1]$ for the function $f_{25}(x,-0.96)$.
Function information is given at $n$ equally spaced points, while
derivative information are given at $n-1$ points in between.
Error has been normalized by
maximum function value.
}
\label{tab:1d_error_tensor_interlace}
\end{table}

\begin{table}
\centering
\begin{tabular}{|c||c|c|c|c|}
\hline
$n$ & $s = 2$ & $s = 3$ & $s = 4$ & $s = 5$ \\
\hline
\hline
11 & 1.42e-01 & 1.03e-01 & 1.04e-01 & 1.02e-01 \\
21 & 5.60e-02 & 2.41e-02 & 1.38e-02 & 8.44e-03 \\
31 & 2.98e-02 & 9.55e-03 & 4.76e-03 & 2.76e-03 \\
41 & 1.88e-02 & 5.48e-03 & 2.15e-03 & 9.72e-04 \\
51 & 1.30e-02 & 3.52e-03 & 1.16e-03 & 4.36e-04 \\
61 & 1.02e-02 & 2.42e-03 & 6.66e-04 & 5.87e-04 \\
\hline
\end{tabular}
\caption{
Interpolation error of MSN 2D Birkhoff interpolation for various $s$
values for the function $f_{25}(x,y)$ defined in Eq.~(\ref{eq:runge_func}).
Function and derivative information is given on $n\times n$ tensor grid.
Error has been normalized by the maximum function value.
All of the computations in this table were performed in
\textbf{single precision}.
}
\label{tab:2d_error_tensor_float32}
\end{table}

\begin{table}
\centering
\begin{tabular}{|c||c|c|c|c|c|c|}
\hline
$n$ & $s = 2$ & $s = 4$ & $s = 6$ & $s = 8$ & $s = 10$ & $s = 12$ \\
\hline
\hline
11 & 1.42e-01 & 1.04e-01 & 1.02e-01 & 3.53e-01 & 9.02e-01 & 1.88e+00 \\
21 & 5.60e-02 & 1.38e-02 & 8.24e-03 & 5.30e-02 & 2.40e-01 & 7.44e-01 \\
31 & 2.98e-02 & 4.76e-03 & 1.72e-03 & 1.33e-03 & 7.60e-03 & 3.75e-02 \\
41 & 1.88e-02 & 2.16e-03 & 4.85e-04 & 1.69e-04 & 1.73e-03 & 1.38e-02 \\
51 & 1.30e-02 & 1.15e-03 & 1.84e-04 & 4.61e-05 & 2.23e-05 & 8.23e-05 \\
61 & 1.02e-02 & 6.81e-04 & 8.19e-05 & 1.59e-05 & 5.57e-06 & 5.59e-05 \\
\hline
\end{tabular}
\caption{
Interpolation error of MSN 2D Birkhoff interpolation for various $s$
values for the function $f_{25}(x,y)$ defined in Eq.~(\ref{eq:runge_func}).
Function and derivative information is given on $n\times n$ tensor grid.
Error has been normalized by the maximum function value.
}
\label{tab:2d_error_tensor}
\end{table}

\begin{table}
\centering
\begin{tabular}{|c||c|c|c|c|c|c|}
\hline
$n$ & $s = 2$ & $s = 4$ & $s = 6$ & $s = 8$ & $s = 10$ & $s = 12$ \\
\hline
\hline
11 & 1.42e-01 & 2.64e-01 & 4.32e-01 & 1.07e+00 & 2.89e+00 & 4.78e+00 \\
21 & 7.57e-02 & 4.13e-02 & 2.03e-02 & 1.67e-01 & 1.13e+00 & 4.41e+00 \\
31 & 3.11e-02 & 9.64e-03 & 1.46e-02 & 2.69e-02 & 1.09e-01 & 3.76e-01 \\
41 & 3.18e-02 & 4.12e-03 & 2.51e-03 & 6.09e-03 & 9.58e-03 & 5.68e-02 \\
51 & 2.05e-02 & 1.54e-03 & 8.55e-04 & 4.63e-04 & 1.32e-03 & 5.02e-03 \\
61 & 9.86e-03 & 9.09e-04 & 1.19e-04 & 4.48e-04 & 4.94e-04 & 4.11e-04 \\
\hline
\end{tabular}
\caption{
Interpolation error of MSN 2D Birkhoff interpolation for various $s$
values for the function $f_{25}(x,y)$ defined in Eq.~(\ref{eq:runge_func}).
Function and derivative information is given on $n\times n$ tensor grid
that has been intersected with an annulus (inner radius $0.5$ and
outer radius $1$).
Error has been normalized by the maximum function value.
}
\label{tab:2d_error_tensor_annulus}
\end{table}

\begin{table}
\centering
\begin{tabular}{|c||c|c|c|c|c|c|}
\hline
$(m,n)$ & $s = 2$ & $s = 4$ & $s = 6$ & $s = 8$ & $s = 10$ & $s = 12$ \\
\hline
\hline
(5,32) & 9.15e-02 & 4.20e-02 & 4.91e-02 & 1.84e-01 & 4.39e-01 & 7.84e-01 \\
(7,48) & 1.74e-02 & 8.52e-03 & 7.92e-03 & 1.97e-02 & 6.31e-02 & 4.05e-01 \\
(9,64) & 7.36e-03 & 2.50e-03 & 1.71e-03 & 4.99e-03 & 1.55e-02 & 6.31e-02 \\
(11,80) & 8.56e-03 & 8.26e-04 & 5.07e-04 & 7.21e-04 & 1.94e-03 & 1.36e-02 \\
(13,96) & 5.62e-03 & 2.85e-04 & 1.73e-04 & 1.73e-04 & 3.70e-04 & 2.22e-03 \\
(15,112) & 2.94e-03 & 1.15e-04 & 6.62e-05 & 3.63e-05 & 5.70e-05 & 3.45e-04 \\
\hline
\end{tabular}
\caption{
Interpolation error of MSN 2D Birkhoff interpolation for various $s$
values for the function $f_{25}(x,y)$ defined in Eq.~(\ref{eq:runge_func}).
Function and derivative information is given on an $(m,n)$ annular grid:
inner radius $0.5$, outer radius $1$, $m$ equally spaced points in the
radial direction, and $n$ equally spaced points in the angular direction.
Error has been normalized by the maximum function value.
}
\label{tab:2d_error_annular_grid}
\end{table}

\begin{table}
\centering
\begin{tabular}{|c||c|c|c|c|}
\hline
$n$ & $s = 2$ & $s = 3$ & $s = 4$ & $s = 5$ \\
\hline
\hline
11 & 1.20e-01 & 8.09e-02 & 6.78e-02 & 1.19e-01 \\
21 & 3.70e-02 & 1.10e-02 & 6.36e-03 & 4.88e-03 \\
31 & 1.88e-02 & 3.00e-03 & 1.02e-03 & 6.21e-04 \\
41 & 1.14e-02 & 1.63e-03 & 4.58e-04 & 1.59e-04 \\
51 & 7.89e-03 & 1.02e-03 & 2.49e-04 & 7.35e-05 \\
61 & 6.04e-03 & 6.87e-04 & 1.48e-04 & 4.12e-05 \\
\hline
\end{tabular}
\caption{
Interpolation error of MSN 2D Birkhoff interpolation
for various $s$ values on $[-1,1]^{2}$ for the function $f_{25}(x,y)$.
Function information is given at $n\times n$ tensor grid, while
derivative information are given at $(n-1)\times(n-1)$ tensor grid in between.
Error has been normalized by
maximum function value.
All of the computations in this table were performed in
\textbf{single precision}.
}
\label{tab:2d_error_tensor_interlace_parallel_float32}
\end{table}

\begin{table}
\centering
\begin{tabular}{|c||c|c|c|c|}
\hline
$n$ & $s = 2$ & $s = 3$ & $s = 4$ & $s = 5$ \\
\hline
\hline
11 & 1.20e-01 & 8.09e-02 & 6.78e-02 & 1.19e-01 \\
21 & 3.70e-02 & 1.10e-02 & 6.36e-03 & 4.88e-03 \\
31 & 1.88e-02 & 3.00e-03 & 1.03e-03 & 6.24e-04 \\
41 & 1.14e-02 & 1.64e-03 & 4.60e-04 & 1.60e-04 \\
51 & 7.89e-03 & 1.02e-03 & 2.47e-04 & 6.73e-05 \\
61 & 6.04e-03 & 6.88e-04 & 1.46e-04 & 3.37e-05 \\
\hline
\end{tabular}
\caption{
Interpolation error of MSN 2D Birkhoff interpolation
for various $s$ values on $[-1,1]^{2}$ for the function $f_{25}(x,y)$.
Function information is given at $n\times n$ tensor grid, while
derivative information are given at $(n-1)\times(n-1)$ tensor grid in between.
\textbf{Directional derivatives are not parallel to the coordinate
axes} and are in the directions
$\left(\frac{1}{\sqrt{2}},\frac{1}{\sqrt{2}}\right)$
and $\left(\frac{1}{\sqrt{2}},-\frac{1}{\sqrt{2}}\right)$.
Error has been normalized by
maximum function value.
All of the computations in this table were performed in
\textbf{single precision}.
}
\label{tab:2d_error_tensor_interlace_not_float32}
\end{table}

\begin{table}
\centering
\begin{tabular}{|c||c|c|c|c|c|c|}
\hline
$n$ & $s = 2$ & $s = 4$ & $s = 6$ & $s = 8$ & $s = 10$ & $s = 12$ \\
\hline
\hline
11 & 1.20e-01 & 6.78e-02 & 3.45e-01 & 1.63e+00 & 4.58e+00 & 1.37e+01 \\
21 & 3.70e-02 & 6.36e-03 & 1.53e-02 & 7.04e-02 & 2.11e-01 & 1.06e+00 \\
31 & 1.88e-02 & 1.03e-03 & 7.85e-04 & 2.11e-03 & 2.82e-02 & 2.71e-01 \\
41 & 1.14e-02 & 4.60e-04 & 6.05e-05 & 8.40e-05 & 3.29e-03 & 4.06e-02 \\
51 & 7.89e-03 & 2.49e-04 & 2.30e-05 & 7.96e-06 & 2.35e-04 & 2.92e-03 \\
61 & 6.04e-03 & 1.47e-04 & 1.12e-05 & 1.68e-06 & 1.39e-05 & 6.29e-05 \\
\hline
\end{tabular}
\caption{
Interpolation error of MSN 2D Birkhoff interpolation
for various $s$ values on $[-1,1]^{2}$ for the function $f_{25}(x,y)$.
Function information is given at $n\times n$ tensor grid, while
derivative information are given at $(n-1)\times(n-1)$ tensor grid in between.
Error has been normalized by
maximum function value.
}
\label{tab:2d_error_tensor_interlace_parallel}
\end{table}

\begin{table}
\centering
\begin{tabular}{|c||c|c|c|c|c|c|}
\hline
$n$ & $s = 2$ & $s = 4$ & $s = 6$ & $s = 8$ & $s = 10$ & $s = 12$ \\
\hline
\hline
11 & 1.20e-01 & 6.78e-02 & 3.45e-01 & 1.63e+00 & 4.58e+00 & 1.37e+01 \\
21 & 3.70e-02 & 6.36e-03 & 1.53e-02 & 7.04e-02 & 2.11e-01 & 1.06e+00 \\
31 & 1.88e-02 & 1.03e-03 & 7.85e-04 & 2.11e-03 & 2.82e-02 & 2.71e-01 \\
41 & 1.14e-02 & 4.60e-04 & 6.05e-05 & 8.40e-05 & 3.29e-03 & 4.06e-02 \\
51 & 7.89e-03 & 2.49e-04 & 2.30e-05 & 7.96e-06 & 2.35e-04 & 2.92e-03 \\
61 & 6.04e-03 & 1.47e-04 & 1.12e-05 & 1.68e-06 & 1.39e-05 & 6.30e-05 \\
\hline
\end{tabular}
\caption{
Interpolation error of MSN 2D Birkhoff interpolation
for various $s$ values on $[-1,1]^{2}$ for the function $f_{25}(x,y)$.
Function information is given at $n\times n$ tensor grid, while
derivative information are given at $(n-1)\times(n-1)$ tensor grid in between.
\textbf{Directional derivatives are not parallel to the coordinate
axes} and are in the directions
$\left(\frac{1}{\sqrt{2}},\frac{1}{\sqrt{2}}\right)$
and $\left(\frac{1}{\sqrt{2}},-\frac{1}{\sqrt{2}}\right)$.
Error has been normalized by
maximum function value.
}
\label{tab:2d_error_tensor_interlace_not}
\end{table}

\begin{table}
\centering
\begin{tabular}{|c||c|c|c|c|}
\hline
$d$ & $s = 2$ & $s = 3$ & $s = 4$ & $s = 5$ \\
\hline
\hline
11 & 2.99e-01 & 2.94e-01 & 3.53e-01 & 4.05e-01 \\
21 & 1.16e-01 & 7.37e-02 & 6.41e-02 & 6.13e-02 \\
31 & 3.63e-02 & 1.38e-02 & 9.05e-03 & 7.65e-03 \\
41 & 1.62e-02 & 3.83e-03 & 1.79e-03 & 1.27e-03 \\
51 & 7.85e-03 & 1.24e-03 & 4.14e-04 & 2.59e-04 \\
\hline
\end{tabular}
\caption{
Interpolation error of MSN 2D Birkhoff interpolation
for various $s$ values on $\SS^{2}$ for the function
$g(x,y,z)$ defined in Eq.~(\ref{eq:runge_func_sphere}).
Function and derivative information are given on a scattered grid
with minimum separation approximately $\pi/d$.
Error has been normalized by the maximum function value.
All of the computations in this table were performed in
\textbf{single precision}.
}
\label{tab:2d_error_sphere_float32}
\end{table}

\begin{table}
\centering
\begin{tabular}{|c||c|c|c|c|c|c|}
\hline
$d$ & $s = 2$ & $s = 4$ & $s = 6$ & $s = 8$ & $s = 10$ & $s = 12$ \\
\hline
\hline
11 & 2.99e-01 & 3.53e-01 & 4.36e-01 & 4.71e-01 & 4.90e-01 & 5.03e-01 \\
21 & 1.16e-01 & 6.41e-02 & 6.04e-02 & 6.07e-02 & 6.14e-02 & 6.19e-02 \\
31 & 3.63e-02 & 9.05e-03 & 7.20e-03 & 7.68e-03 & 8.35e-03 & 9.23e-03 \\
41 & 1.62e-02 & 1.79e-03 & 1.06e-03 & 9.16e-04 & 8.90e-04 & 8.68e-04 \\
51 & 7.85e-03 & 4.12e-04 & 1.89e-04 & 1.64e-04 & 1.59e-04 & 1.59e-04 \\
\hline
\end{tabular}
\caption{
Interpolation error of MSN 2D Birkhoff interpolation
for various $s$ values on $\SS^{2}$ for the function
$g(x,y,z)$ defined in Eq.~(\ref{eq:runge_func_sphere}).
Function and derivative information are given on a scattered grid
with minimum separation approximately $\pi/d$.
Error has been normalized by the maximum function value.
}
\label{tab:2d_error_sphere}
\end{table}

\begin{table}
\centering
\begin{tabular}{|c||c|c|c|c|c|c|}
\hline
$d$ & $s = 2$ & $s = 4$ & $s = 6$ & $s = 8$ & $s = 10$ & $s = 12$ \\
\hline
\hline
15 & 6.80e-01 & 9.08e-01 & 1.20e+00 & 1.54e+00 & 1.78e+00 & 1.89e+00 \\
39 & 2.35e-01 & 2.52e-01 & 2.18e-01 & 1.60e-01 & 2.04e-01 & 3.65e-01 \\
63 & 9.12e-02 & 4.83e-02 & 4.21e-02 & 3.83e-02 & 3.50e-02 & 2.80e-02 \\
87 & 1.92e-02 & 7.11e-03 & 6.49e-03 & 5.89e-03 & 4.75e-03 & 3.57e-03 \\
111 & 1.71e-02 & 1.89e-03 & 1.38e-03 & 1.03e-03 & 6.48e-04 & 5.94e-04 \\
\hline
\end{tabular}
\caption{
Interpolation error of MSN 2D Birkhoff interpolation
for various $s$ values on $\SS^{2}$ for the function
$g(x,y,z)$ defined in Eq.~(\ref{eq:runge_func_sphere}).
Function and derivative information are given on a scattered grid
with minimum separation approximately $\pi/d$ with the restriction
that the polar angle $\theta\in[0,\frac{\pi}{3}]\cup[\frac{2\pi}{3},\pi]$.
Error has been normalized by the maximum function value.
}
\label{tab:2d_error_sphere_north_south}
\end{table}

\begin{table}
\centering
\begin{tabular}{|c||c|c|c|c|}
\hline
$n$ & $s = 2$ & $s = 3$ & $s = 4$ & $s = 5$ \\
\hline
\hline
11 & 6.23e-02 & 3.15e-02 & 2.55e-02 & 2.03e-02 \\
21 & 2.20e-02 & 6.85e-03 & 2.79e-03 & 1.36e-03 \\
31 & 1.12e-02 & 2.78e-03 & 8.05e-04 & 2.96e-04 \\
41 & 7.52e-03 & 1.50e-03 & 3.35e-04 & 1.43e-04 \\
51 & 5.52e-03 & 9.08e-04 & 1.70e-04 & 8.42e-05 \\
61 & 4.24e-03 & 5.95e-04 & 6.60e-04 & 8.80e-04 \\
\hline
\end{tabular}
\caption{
Interpolation error of MSN 2D Birkhoff interpolation for various $s$
values for the function $f_{9}(x,y)$ defined in Eq.~(\ref{eq:runge_func}).
Function and derivative information is given on $n\times n$ tensor grid.
Error has been normalized by the maximum function value.
All of the computations in this table were performed in
\textbf{single precision}.
}
\label{tab:2d_error_tensor_float32_easy}
\end{table}

\begin{table}
\centering
\begin{tabular}{|c||c|c|c|c|c|c|}
\hline
$n$ & $s = 2$ & $s = 4$ & $s = 6$ & $s = 8$ & $s = 10$ & $s = 12$ \\
\hline
\hline
11 & 6.23e-02 & 2.55e-02 & 2.66e-02 & 8.11e-02 & 1.72e-01 & 3.00e-01 \\
21 & 2.20e-02 & 2.79e-03 & 7.54e-04 & 1.48e-03 & 6.72e-03 & 2.05e-02 \\
31 & 1.12e-02 & 8.06e-04 & 1.19e-04 & 3.49e-05 & 2.97e-05 & 1.81e-04 \\
41 & 7.51e-03 & 3.30e-04 & 2.91e-05 & 5.55e-06 & 1.73e-06 & 5.67e-06 \\
51 & 5.52e-03 & 1.63e-04 & 9.66e-06 & 1.37e-06 & 2.83e-07 & 9.92e-08 \\
61 & 4.24e-03 & 9.00e-05 & 3.89e-06 & 4.34e-07 & 6.54e-08 & 1.75e-08 \\
\hline
\end{tabular}
\caption{
Interpolation error of MSN 2D Birkhoff interpolation for various $s$
values for the function $f_{9}(x,y)$ defined in Eq.~(\ref{eq:runge_func}).
Function and derivative information is given on $n\times n$ tensor grid.
Error has been normalized by the maximum function value.
}
\label{tab:2d_error_tensor_easy}
\end{table}


%

\bhag{Conclusions}
We proved a general Birkhoff interpolation result: with
minimal restrictions it is possible
to interpolate function and derivative information with diffusion
polynomials of degree $N^{*} \sim \eta^{-1}$, where $\eta$ is the
minimum separation distance between points.
This extends previous work related to function interpolation~\cite{bdint}
and is needed for numerical approximations arising in the solution
of partial differential equations~\cite{chandra2015minimum}.
One and two dimensional
numerical examples were presented to demonstrate the abilities of this
method with indifference to interpolation location.

\bibliographystyle{plain}
\bibliography{hrushikesh_mod.bib}
\end{document}